\documentclass[11pt]{amsart}
\usepackage{CJK}
\usepackage[centertags]{amsmath}
\usepackage{times,fullpage}
\usepackage{latexsym,amscd,amssymb}
\usepackage{tikz,tikz-cd}
\usepackage{dutchcal}
\usepackage[all]{xy}
 \input xy
  \xyoption{all}

\usepackage{amsfonts}
\usepackage{amssymb}
\usepackage{amsthm}
\usepackage{newlfont}
\usepackage[colorlinks,linkcolor=blue, pdfstartview=FitH]{hyperref}
\usepackage{mathrsfs}
\usepackage{enumerate}
\numberwithin{equation}{section}

 \theoremstyle{plain}
\newtheorem{thm}{Theorem}[section]
\newtheorem{theorem}[thm]{Theorem}
\newtheorem{lemma}[thm]{Lemma}
\newtheorem{corollary}[thm]{Corollary}
\newtheorem{proposition}[thm]{Proposition}

\theoremstyle{definition}

\newtheorem{question}[thm]{Question}
\newtheorem{conjecture}[thm]{Conjecture}

\newtheorem{remark}[thm]{Remark}

\newtheorem{definition}[thm]{Definition}
\newtheorem{claim}[thm]{Claim}

\newtheorem{example}[thm]{Example}

\newtheorem{defn-thm}[thm]{Definition-Theorem}

\newcommand{\sF}{{\mathcal F}}

\newcommand{\sO}{{\mathcal O}}

\newcommand{\E}{{\mathbb E}}

\newcommand{\N}{{\mathbb N}}

\newcommand{\Q}{{\mathbb Q}}

\newcommand{\btheorem}{\begin{theorem}}
\newcommand{\etheorem}{\end{theorem}}
\newcommand{\bquestion}{\begin{question}}
\newcommand{\equestion}{\end{question}}
\newcommand{\bconjecture}{\begin{conjecture}}
\newcommand{\econjecture}{\end{conjecture}}
\newcommand{\bclaim}{\begin{claim}}
\newcommand{\eclaim}{\end{claim}}
\newcommand{\bproposition}{\begin{proposition}}
\newcommand{\eproposition}{\end{proposition}}
\newcommand{\bdefinition}{\begin{definition}}
\newcommand{\edefinition}{\end{definition}}
\newcommand{\bcorollary}{\begin{corollary}}
\newcommand{\ecorollary}{\end{corollary}}
\newcommand{\bproof}{\begin{proof}}
\newcommand{\eproof}{\end{proof}}
\newcommand{\bremark}{\begin{remark}}
\newcommand{\eremark}{\end{remark}}
\newcommand{\eexample}{\end{example}}
\newcommand{\bexample}{\begin{example}}
\newcommand{\elemma}{\end{lemma}}
\newcommand{\blemma}{\begin{lemma}}

\renewcommand{\phi}{\varphi}

\newcommand{\ee}{\end{eqnarray*}}
\newcommand{\be}{\begin{eqnarray*}}

\newcommand{\beq}{\begin{equation}}
\newcommand{\eeq}{\end{equation}}

\newcommand{\bd}{\begin{enumerate}[(1)]}
\newcommand{\ed}{\end{enumerate}}

\renewcommand{\hat}{\widehat}
\renewcommand{\tilde}{\widetilde}

\renewcommand{\bf}{\textbf}

\renewcommand{\rm}{\textrm}

\newcommand{\ratmap}{\dashrightarrow}

\begin{document}

\title{Note on the 3-dimensional log canonical abundance in characteristic $>3$}
\author{Zheng Xu}

\address{Academy of Mathematics and Systems Science, Chinese Academy of Sciences,
Beijing, China}
\email{zxu@amss.ac.cn}

\date{}

\maketitle

\begin{abstract}
In this paper, we prove the non-vanishing and some special cases of the abundance for log canonical threefold pairs over an algebraically closed field $k$ of characteristic $p > 3$. More  precisely, we prove that if $(X,B)$ be a projective log canonical threefold pair over $k$ and $K_{X}+B$ is pseudo-effective, then $\kappa(K_{X}+B)\geq 0$, and if $K_{X}+B$ is nef and $\kappa(K_{X}+B)\geq 1$, then $K_{X}+B$ is semi-ample. 

As applications, we show that the log canonical rings of projective log canonical threefold pairs over $k$ are finitely generated and the abundance holds when the nef dimension $n(K_{X}+B)\leq 2$ or when the Albanese map  $a_{X}:X\to \mathrm{Alb}(X)$ is non-trivial. Moreover, we prove that the abundance for klt threefold pairs over $k$ implies the abundance for log canonical threefold pairs over $k$.
\end{abstract}

{\tableofcontents \setcounter{page}{0}\pagenumbering{roman}}

\setcounter{page}{1}\pagenumbering{arabic}

\section{Introduction}

Over the last decade, the Minimal Model Program (MMP) for threefolds over a field of characteristic $>3$ has been largely established.
First, Hacon and Xu proved the existence of minimal models for terminal threefolds over an algebraically closed field $k$ of characteristic $ > 5$ (\cite{hacon2015three}). Then Cascini, Tanaka and Xu proved that arbitrary terminal threefold over $k$ is birational to either a minimal model or a Mori fibre space (\cite{cascini2015base}). Base on it, Birkar and Waldron established
the MMP for klt threefolds over $k$ (\cite{birkar2016existence,birkar2017existence}). Moreover, there are some generalizations of it in various directions. For example, see \cite{hashizume2020minimal,waldron2018lmmp} for its generalization to log canonical (lc) pairs, \cite{gongyo2019rational,hacon2022minimal,hacon2022relative} for its generalization to low characteristics, \cite{das2022log} for its generalization to imperfect base fields, and \cite{bhatt2023globally} for its analog in mixed characteristics.

Now we can run MMPs for lc threefold pairs over a perfect field of characteristic $>3$ (see Theorem \ref{lcmmp}). Hence a central problem remaining is the following conjecture.\\

\noindent \bf{Abundance conjecture.} 
Let $(X,B)$ be a projective lc threefold pair over a perfect field $k$ of characteristic $>3$. If $K_{X}+B$ is nef,
then it is semi-ample.

\bremark The abundance conjecture for lc surface pairs over any field of positive characteristic is proved in \cite{tanaka2020abundance}, and for slc surface pairs over any field of positive characteristic it is proved in \cite{posva2023abundance}.
\eremark

\bremark (From a perfect field to its algebraic closure) Many properties of singularities and positivity, e.g. klt, lc, semi-ampleness and Iitaka dimensions, are preserved under the base change from a perfect field to its algebraic closure (see \cite[Remark 2.7]{gongyo2019rational} for example). In this paper, we sometimes do such base changes and assume that we work over algebraically closed fields. However, some conditions need that the base field is algebraically closed, e.g. conditions about nef dimensions (see Subsection \ref{nefreductionmap} for definition) and Albanese maps.
\eremark

When $K_{X}+B$ is big, Birkar and Waldron proved it in characteristic $>5$ (\cite{birkar2017existence,waldron2018lmmp}), then Hacon and Witaszek proved it in characteristic $5$ (\cite{hacon2022minimal}). 
When $(X,B)$ is klt and the characteristic of $k$ is greater than $5$, 
Waldron proved it in the case of $\kappa(X,K_{X}+B)= 2$ (\cite{waldron2017finite}), Das, Waldron and Zhang proved it in the case of $\kappa(X,K_{X}+B)= 1$ (\cite{das2019abundance,zhang2019abundance}), Witaszek proved it in the case when the nef dimension $n(X,K_{X}+B)\leq 2$ (\cite{witaszek2021canonical}), and Zhang proved it in the case when the Albanese map $a_{X}: X\to \mathrm{Alb}(X)$ is non-trivial (\cite{zhang2020abundance}).
In conclusion, the abundance holds when $(X,B)$ is klt, the characteristic of $k$ is greater than $5$ and one of the following conditions holds:

\noindent (1) $\kappa(X,K_{X}+B)\geq 1$,
 
\noindent (2) the nef dimension $n(X,K_{X}+B)\leq 2$,

\noindent (3) the Albanese map $a_{X}: X\to \mathrm{Alb}(X)$ is non-trivial.

The above works on the abundance for klt pairs in characteristic $>5$ can be generalized to the case when the characteristic is greater than $3$ by some careful modifications (see Section 3).
Then it is natural to ask the following question.

\bquestion How can we generalize a result on the abundance for klt threefold pairs to lc threefold pairs?
\equestion

In characteristic $0$, this is done in \cite{keel1994log}. However, the approach there needs vanishing theorems and the termination of flips for threefolds. The vanishing theorems may fail in positive characteristic and the termination of flips for threefolds is unknown in positive characteristic for lack of a good understanding of terminal threefold singularities in positive characteristic.
In this paper, we propose a new method to solve Question 1.1 and  generalize most of results on the abundance for klt pairs in characteristic $>5$ to lc pairs in characteristic $>3$.
We first prove the nonvanishing theorem for lc threefold pairs over a perfect field $k$ of characteristic $>3$.

\btheorem (Theorem \ref{nonvan}) Let $(X,B)$ be a projective lc threefold pair over a perfect field $k$ of characteristic $>3$. If $K_{X}+ B$ is pseudo-effective, then $\kappa(X,K_{X}+B)\geq 0$.
\etheorem

As a corollary, we have the following result on termination of flips.

\btheorem (Theorem \ref{termination}) Let $(X,B)$ be a projective lc threefold pair defined over a perfect field $k$ of characteristic $p > 3$ such that $K_{X}+B$ is pseudo-effective. Then every sequence of $(K_{X}+B)$-flips terminates. In particular, any $(K_{X}+B)$-MMP terminates with a minimal model.
\etheorem

Secondly, we prove the following result which is the main technical result of this paper.

\btheorem (Theorem \ref{mainthm}) Let $(X,B)$ be a projective lc threefold pair over an algebraically closed field  $k$ of characteristic $>3$. If $K_{X}+B$ is nef and $\kappa(X,K_{X}+B)\geq 1$, then $K_{X}+B$ is semi-ample.
\etheorem

Combined with the results on klt pairs, we deduce the following statements.

\btheorem (Theorem \ref{finitegeneration}) Let $(X, B)$ be a projective lc threefold pair over an algebraically closed field  $k$ of characteristic $>3$. Then the log canonical ring
$$R(K_{X}+B)=\oplus_{m=0}^{\infty} H^{0}(\lfloor m(K_{X}+B)\rfloor )$$
is finitely generated.
\etheorem

\btheorem (Theorem \ref{nefdim}) Let $(X,B)$ be a projective lc threefold pair over an algebraically closed field  $k$ of characteristic $>3$. If $K_{X}+B$ is nef and the nef dimension $n(X,K_{X}+B)\leq 2$, then $K_{X}+B$ is semi-ample.
\etheorem

\btheorem (Theorem \ref{Alb}) 
 Let $(X,B)$ be a projective lc threefold pair over an algebraically closed field $k$ of characteristic $>3$.  If $K_{X}+B$ is nef and $\mathrm{dim}\ \mathrm{Alb}(X)\neq 0$, then $K_{X}+B$ is semi-ample.
\etheorem

It turns out that the following result follows from Theorem 1.4 and Theorem 1.6.

\btheorem (Theorem \ref{kltimplylc}) Let $k$ be an algebraically closed field of characteristic $>3$. 
Assume we have

\noindent (1) abundance for terminal threefolds over $k$ holds, and

\noindent (2) any effective nef divisor $D$ on any klt Calabi-Yau threefold pair $(Y,\Delta)$ ($(Y,\Delta)$ is klt and $K_{Y}+\Delta\sim_{\Q}0$) over $k$ is semi-ample.

\noindent Then the abundance conjecture for threefold pairs over $k$ holds.
In particular, the abundance conjecture for klt threefold pairs over $k$ implies the abundance conjecture for lc threefold pairs over $k$.
\etheorem

\noindent \bf{Outline of the proof of Theorem 1.6.}

For simplicity, we assume that $k$ is an uncountable algebraically closed field of characteristic $>3$ (the uncountability is used for defining the nef reduction map).
We first prove the nonvanishing theorem for projective lc threefold pairs over $k$ (see Theorem \ref{nonvan}) as follows.
By Theorem \ref{dltmodification}, after replacing, we can assume that $(X,B)$ is $\Q$-factorial and dlt, and moreover $X$ is terminal. Then we run a $K_{X}$-MMP which is $(K_{X}+B)$-trivial by  Definition \ref{MMPKtrivial}. It terminates by Lemma \ref{terminaltermination}. If we get a minimal model, then we can use the nonvanishing for klt pairs (see Theorem \ref{kltnonvan}) to prove the assertion. Otherwise, we get a Mori fibre space. It implies that the nef dimension $n(K_{X}+B)\leq 2$. We can use Witaszek's weak canonical bundle formula to handle the case of $n(K_{X}+B)= 2$. The case of $n(K_{X}+B)= 1$ is trivial by descenting $K_{X}+B$ along the nef reduction map of $K_{X}+B$. Finally, we need to handle the case of $n(K_{X}+B)= 0$. In this case, $K_{X}+B$ is numerically trivial. Then the semi-ampleness of $K_{X}+B$ preserves under any step of MMPs. By Theorem \ref{lcmmp}, we can run a $(K_{X}+B-\lfloor B\rfloor)$-MMP which terminates. It terminates with a Mori fibre space and then we can descent $K_{X}+B$ along the Mori fibre space to prove its semi-ampleness. In conclusion, the nonvanishing holds.
As a corollary, we have the termination of flips for pseudo-effective lc threefold pairs over $k$ (see Theorem \ref{termination}).

Now let $(X,B)$ be a projective lc threefold pair over $k$ such that $K_{X}+B$ is nef.
We assume $\kappa(K_{X}+B)=2$, which is the most difficult case. Then $K_{X}+B$ is endowed with a map $h:X\to Z$ to a normal proper algebraic space of dimension $2$ by Lemma \ref{diagram}.
We replace $(X,B)$ by a $\Q$-factorial dlt modification by Theorem \ref{dltmodification}.
Then one of the following cases holds:\\

\noindent Case \uppercase\expandafter{\romannumeral1}: $K_{X}+B-\varepsilon\lfloor B\rfloor$ is not pseudo-effective for any rational $\varepsilon>0$,\\

\noindent Case \uppercase\expandafter{\romannumeral2}: $K_{X}+B-\varepsilon\lfloor B\rfloor$ is pseudo-effective for any sufficiently small rational $\varepsilon>0$.\\

In Case \uppercase\expandafter{\romannumeral1}, we first prove that $\lfloor B\rfloor$ must dominate $Z$ (see Proposition \ref{terminationwhendegenerated}).
Then we deduce the semi-ampleness of $K_{X}+B$ by adjunction (see Proposition \ref{dominatedcase}).

In Case \uppercase\expandafter{\romannumeral2}, we first modify the pair $(X,B)$ by running several MMP which are $(K_{X}+B)$-trivial (see Definition \ref{MMPKtrivial}) so that all $h$-exceptional prime divisors are connected components of $\lfloor B\rfloor$.
Then after further modification we can construct an equidimensional fibration
$h_{\varepsilon}:X\to Z_{\varepsilon}$ to a normal projective surface.
Finally, we descend $K_{X}+B$ to $Z_{\varepsilon}$ and prove its semi-ampleness (see Proposition \ref{Case2}).
$\hfill\square$\\

\noindent\bf{Notation and conventions.}

$\bullet$ We say that $X$ is a variety if it is an integral and separated scheme which is  of finite type over a field $k$.

$\bullet$ We say that a morphism $f: X\to Y$ is a contraction if $X$ and $Y$ are normal
algebraic spaces (we refer to \cite{artin1971algebraic} for definition and basic properties of algebraic spaces), $f_{\ast}\sO_{X}=\sO_{Y}$, and $f$ is proper.

$\bullet$ We say that a morphism $f : X\to Y$ of algebraic spaces is equidimensional if all fibres $X_{y}$ of $f$ are of the same dimension for $y\in Y$.

$\bullet$ Let $f : X\to Y$ be a surjective morphism of integral algebraic spaces. 
We say that a $\Q$-divisor $D$ on $X$ is $f$-exceptional if $\mathrm{dim}(f(\mathrm{Supp}\ D))<\mathrm{dim}\ Y-1$.

$\bullet$ We call a divisor $D \subseteq X$ vertical with respect to a contraction $f$ if $f|_{D}$ is not dominant.

$\bullet$ We call $(X,B)$ a pair if $X$ is a normal variety and $B$ is an effective $\Q$-divisor on $X$ such that $K_{X}+B$ is $\Q$-Cartier. For more notions in the theory of MMP 
such as klt (dlt, lc) pairs, filps, divisorial contractions and so on, we refer to \cite{kollar1998birational}.

$\bullet$ Let $X$ be a normal projective variety over a field $k$ and $D$ be a $\Q$-Cartier $\Q$-divisor on $X$. If $|mD|=\emptyset$ for all $m>0$, we define the Kodaira dimension $\kappa(X,D)=-\infty$. Otherwise, let $\Phi:X\dashrightarrow Z$ be the Iitaka map (we refer to \cite[2.1.C]{lazarsfeld2017positivity}) of $D$ and we define the Kodaira dimension $\kappa(X,D)$ to be the dimension of the image of $\Phi$. Sometimes we write $\kappa(D)$ for $\kappa(X,D)$. We denote $\kappa(X,K_{X})$ by $\kappa(X)$. And for a projective variety $Y$ over a field $k$ admitting a smooth model $\tilde{Y}$, we define $\kappa(Y):=\kappa(\tilde{Y})$.

$\bullet$ Let $X$ be a normal projective variety of dimension $n$ over a field $k$ and $D$ be a   nef $\Q$-Cartier $\Q$-divisor on $X$.
Then we can define
$$\nu(D):=\mathrm{max}\{k\in\N | D^{k}\cdot A^{n-k}>0\ \rm{for}\ \rm{an}\ \rm{ample}\ \rm{divisor}\ A\ \rm{on}\ X\}.$$
\\

\noindent\bf{Acknowledgements.}
I would like to express my gratitude to my advisor Wenhao Ou for his help,
encouragement, and support. Further, I would like to thank Jakub Witaszek for answering my question on his paper \cite{witaszek2021canonical}. Finally, I would like to thank Lei Zhang for his encouragement and helpful advice.

\section{Preliminaries}

In this section we recall some basic results.

\subsection{Keel's results on semi-ampleness}

In this subsection, we survey Keel's work on basepoint free theorem for nef and big $\Q$-Cartier $\Q$-divisors in positive characteristic (see \cite{keel1999basepoint}). It is proved that to show the semi-ampleness of a nef and big $\Q$-Cartier $\Q$-divisor $L$ on a projective variety $X$, it suffices to show the semi-ampleness of $D$ on $\E(L)$, which is a closed subset of $X$ defined below.

\bdefinition  Let $L$ be a nef $\Q$-Cartier $\Q$-divisor on a projective scheme $X$ over a
 field. An irreducible subvariety $Z \subset X$ is called exceptional for $L$ if $L|_{Z}$  is not big, i.e. if $L^{\mathrm{dim}\ Z}\cdot Z = 0$. The exceptional locus of $L$, denoted by $\E(L)$, is the closure of the union of all exceptional subvarieties.
\edefinition
\bremark 
 $\E(L)$ is actually the union of finitely many exceptional subvarieties by \cite[1.2]{keel1999basepoint}.
\eremark

\bdefinition  A nef $\Q$-Cartier $\Q$-divisor $L$ on a proper scheme $X$ over a field is endowed with a map (EWM) $f: X\to Z$ if $f$ is a proper map  to a proper algebraic space $Z$ such that it contracts a closed subvariety $Y$, i.e. $\mathrm{dim}(f(Y))<\mathrm{dim}(Y)$, if and only if $L|_{Y}$ is not big. We may always assume that such a map has geometrically connected fibres.
\edefinition
\bremark By definition, if $L$ is endowed with a map $f: X\to Z$, then a curve $C\subseteq X$ is contracted by $f$ if and only if $L\cdot C=0$.
Moreover, if $f^{\prime}:X\to Z^{\prime}$ is a contraction which only contracts $L$-numerically trivial curves, then by the rigidity lemma (\cite[\uppercase\expandafter{\romannumeral2}.5.3]{kollar2013rational}) $f$ factors through $f^{\prime}$.
\eremark

\blemma\label{pullbacksemi-ample}  Let $p: Y\to X$ be a proper surjective morphism between reduced algebraic spaces of finite type over a field of positive characteristic. Let $L$ be a $\Q$-Cartier $\Q$-divisor on $X$ such that $p^{\ast}L$ is semi-ample. If $X$ is normal, then $L$ is semi-ample.
\elemma
\bproof This lemma follows from \cite[Lemma 2.10]{keel1999basepoint}.
\eproof

The following theorem is the main result of \cite{keel1999basepoint}.

\btheorem\label{Keelsemi-ample} (\cite[Theorem 0.2]{keel1999basepoint}) Let $L$ be a nef $\Q$-Cartier $\Q$-divisor on a scheme $X$, projective over a field of positive characteristic. Then $L$ is semi-ample (resp. EWM) if and only if $L|_{\E(L)}$ is semi-ample (resp. EWM).
\etheorem

\subsection{Nef reduction map}\label{nefreductionmap}

In this subsection, we recall the notion of nef reduction map.

\bdefinition Let $X$ be a normal projective variety defined over an uncountable field
 and let $L$ be a nef $\Q$-Cartier  $\Q$-divisor. We call a rational map 
$\phi: X \dashrightarrow Z$ a nef reduction map of $L$ if $Z$ is a normal projective variety and there exist open dense subsets $U \subseteq X$, $V\subseteq Z$ such that

\noindent (1)  $\phi|_{U}: U \to Z$ is proper,  its image is $V$  and 
$\phi_{\ast}\sO_{U} = \sO_{V}$ ,

\noindent (2)  $L|_{F}\equiv 0$ for all fibres $F$ of $\phi$ over $V$, and

\noindent (3)  if $x \in X$ is a very general point and $C$ is a curve passing through it, then $C\cdot L = 0$ if and only if $C$ is contracted by $\phi$. 
\edefinition

It is proved that a nef reduction map exists over an uncountable algebraically closed field.

\btheorem\label{nefreductionmap} (\cite[Theorem 2.1]{bauer2002reduction}) A nef reduction map exists for normal projective varieties defined over an uncountable algebraically closed field. Furthermore, it is unique
up to birational equivalence.
\etheorem

For a nef reduction map $\phi: X \dashrightarrow Z$ of $L$, 
the nef dimension of $L$ is defined
to be $\mathrm{dim}\ Z$ and denoted by $n(X,L)$. When the base field is countable and algebraically closed, we can define 
$$n(X,L):=n(X_{K},L_{K})$$ 
by \cite[Proposition 2.16]{witaszek2021canonical}, where $K$ is an uncountable algebraically closed field that contains $k$, and $X_{K},L_{K}$ are the base changes of $X,L$ to $K$.
It satisfies $\kappa(X, L)\leq n(X, L)$. Sometimes we write $n(L)$ for $n(X,L)$.

\blemma\label{EWM} (\cite[Lemma 7.2]{birkar2017existence}) Let $X$ be a normal projective variety of dimension $\leq 3$ over an uncountable algebraically closed field of characteristic $p> 0$. Suppose $L$ is a nef $\Q$-Cartier $\Q$-divisor on $X$ with  $\kappa(L)=n(L)\leq 2$. Then $L$ is $EWM$ to a proper algebraic space $Z$ of dimension equal to $\kappa(L)$.
\elemma

The following lemma is very useful for descending a nef $\Q$-Cartier $\Q$-divisor along a fibration.

\blemma\label{descend} Let $f: X\to Z$ be a projective contraction between normal quasi-projective varieties over a field of characteristic $p>0$ and $L$ a $f$-nef $\Q$-Cartier $\Q$-divisor on $X$ such that
$L|_{F} \sim_{\Q} 0$, where $F$ is the generic fibre of $f$. Assume $\mathrm{dim}\ Z\leq 3$. Then there exists a diagram
\beq\begin{tikzcd}
X^{\prime}\arrow[d,swap,"f^{\prime} "]\arrow[r,"\phi "]&X\arrow[d,"f"] \\
Z^{\prime}\arrow[r,"\psi "] &Z\notag
\end{tikzcd}
\eeq
with $\phi,\psi$ projective birational, and a $\Q$-Cartier $\Q$-divisor $D$ on $Z^{\prime}$ such that $\phi^{\ast}L \sim_{\Q} f^{\prime\ast}D$. Moreover, if $Z$ is $\Q$-factorial and $f$ is equidimensional,
then we can take $X^{\prime} = X$ and $Z^{\prime} = Z$.
\elemma
\bproof It is an adaptation of a result of Kawamata \cite[Proposition 2.1]{kawamata1985pluricanonical}. See \cite[Lemma 3.2]{waldron2017finite} for a proof in this setting.
\eproof

\subsection{Abundance theorem for surfaces}

Abundance for slc surfaces over an arbitrary field of characteristic $>0$ is known.

\btheorem\label{surfaceabundance} (\cite[Theorem 1]{posva2023abundance}) Let $(X,\Delta)$ be a projective slc surface pair over a field of characteristic $>0$. If $K_{X}+\Delta$ is nef, then it is semi-ample.
\etheorem

\subsection{MMP for threefolds in positive characteristic}

In this subsection, we recall the theory of MMP for projective lc threefold pairs over a perfect field of characteristic $p> 3$. Moreover, we define a partial MMP over an algebraically closed field of characteristic $p> 3$ (see Definition \ref{MMPKtrivial}). We will use this construction to study the abundance in Section $5$.
 
\btheorem\label{lcmmp} (\cite[Theorem 1.1]{hashizume2020minimal} and \cite{hacon2022minimal})  Let $(X, B)$ be a lc threefold pair over a perfect field $k$ of characteristic $>3$ and $f : X \to Y$ a projective surjective morphism to a quasi-projective variety. If $K_{X}+B$ is pseudo-effective (resp. not pseudo-effective) over $Y$, then we can run a $(K_{X}+B)$-MMP over $Y$ to get a log minimal model (resp. Mori fibre space) over $Y$ .
\etheorem

We recall the notion of MMP with scaling.
Let $(X,B)$ be a projective lc threefold pair over a perfect field $k$ of characteristic $>3$ and $A> 0$ an $\Q$-Cartier $\Q$-divisor on $X$. Suppose that there is $t_{0}> 0$ such that $(X,B+ t_{0}A)$ is lc and $K_{X}+B+t_{0}A$ is nef. We describe how to run a $(K_{X}+B)$-MMP with scaling of $A$ as follows.

Let $\lambda_{0}=\mathrm{inf} \{t|\ K_{X}+B+tA\ \mathrm{is}\ \mathrm{nef} \}$. Suppose we can find a $(K_{X}+ B)$-negative extremal ray $R_{0}$ which satisfies $(K_{X}+B+\lambda_{0}A)\cdot R_{0}= 0$ (In general, it is possible that there is no such extremal ray). This is the first ray we contract in our MMP. If the contraction
is a Mori fibre contraction, we stop. Otherwise let $X_{1}$ be the result of the
divisorial contraction or flip. Then $K_{X_{1}}+ B_{X_{1}}+\lambda_{0}A_{X_{1}}$
is also nef, where $B_{X_{1}}$ and $A_{X_{1}}$ denote the birational transforms on $X_{1}$ of $B$ and $A$, respectively. We define $\lambda_{1} =\mathrm{inf} \{t|\ K_{X_{1}}+ B_{X_{1}} + tA_{X_{1}}\
\mathrm{is}\ \mathrm{nef} \}$. The next step in our MMP is chosen to be a $(K_{X_{1}} + B_{X_{1}})$-negative extremal ray $R_{1}$ which is $(K_{X_{1}}+B_{X_{1}}+ \lambda_{1} A_{X_{1}})$
-trivial. So long as we can find the appropriate extremal
rays, contractions and flips, we can continue this process.

\bproposition\label{findextremalray} Let $(X,B)$ be a $\Q$-factorial projective lc threefold pair over an algebraically closed field $k$ of characteristic $>3$ and $W$ be an effective $\Q$-divisor such that
$K_{X}+B+W$ is nef. Then either

\noindent (1) there is a $(K_{X}+B)$-negative extremal ray which is $(K_{X}+B+W)$-trivial, or

\noindent (2) $K_{X}+B+(1-\varepsilon)W$ is nef for any sufficiently small rational $\varepsilon>0$.
\bproof It is an adaptation of \cite[Lemma 5.1]{keel1994log}. Note that the proof there only uses the fact that for any $(K_{X}+B)$-negative extremal ray $R$ there is a rational curve $C$ such that $C$ generates $R$ and $-(K_{X}+B)\cdot C\leq 6$, which holds in our setting by \cite[Theorem 1.3]{hashizume2020minimal} and \cite{hacon2022minimal}.
\eproof
\eproposition
\bremark The assumption that $k$ is algebraically closed is used for the fact that for any $(K_{X}+B)$-negative extremal ray $R$ there is a rational curve $C$ such that $C$ generates $R$ and $-(K_{X}+B)\cdot C\leq 6$.
\eremark

\bcorollary\label{MMPwithscaling} Let $(X,B)$ be a $\Q$-factorial projective lc threefold pair over an algebraically closed field $k$ of characteristic $>3$ and $A$ be an effective $\Q$-divisor such that $(X,B+A)$ is lc and $K_{X}+B+A$ is nef.
If $K_{X}+B$ is not nef, then we can run a $(K_{X}+B)$-MMP with scaling of $A$.
\bproof Let $\lambda:= \mathrm{inf} \{t|\ K_{X}+B+tA\ \mathrm{is}\ \mathrm{nef} \}$ be the 
nef threshold. Then the only assertion is that we can find a $(K_{X}+B)$-negative extremal ray $R$ such that $(K_{X}+B+\lambda A)\cdot R=0$. We apply Proposition \ref{findextremalray} by letting $W:=\lambda A$.
\eproof
\ecorollary

In this paper, we will use the following construction.

\bdefinition\label{MMPKtrivial} Let $(X,B)$ be a $\Q$-factorial projective lc threefold pair over an algebraically closed field $k$ of characteristic $>3$ and $A$ be an effective $\Q$-divisor such that $(X,B+A)$ is lc and $K_{X}+B+A$ is nef. We can run a partial $(K_{X}+B)$-MMP with scaling of $A$ as follows.

Let $\lambda_{0}=\mathrm{inf} \{t |\ K_{X}+B+tA\ \mathrm{is}\ \mathrm{nef} \}$. If $\lambda_{0}<1$, then we stop. Otherwise, by Proposition \ref{findextremalray} there exists a $(K_{X}+ B)$-negative extremal ray $R_{0}$ which satisfies $(K_{X}+B+A)\cdot R_{0}= 0$. We contract this extremal ray. If the contraction
is a Mori fibre contraction, we stop. Otherwise let $\mu_{0}: X\dashrightarrow X_{1}$ be the divisorial contraction or flip. Repeat this process for $(X_{1},\mu_{0\ast}B), \mu_{0\ast}A$ and so on.

We call this construction a $(K_{X}+B)$-MMP which is $(K_{X}+B+A)$-trivial.
\edefinition

The following lemma tells us what the output of this construction is if it terminates.

\blemma\label{Output} Let $(X,B)$ be a $\Q$-factorial projective lc threefold pair over an algebraically closed field $k$ of characteristic $>3$ and $A$ be an effective $\Q$-divisor such that $(X,B+A)$ is lc and $K_{X}+B+A$ is nef. 

If a $(K_{X}+B)$-MMP which is $(K_{X}+B+A)$-trivial terminates, then its output is a $\Q$-factorial projective lc pair $(X^{\prime}, B^{\prime}+A^{\prime})$, and either

\noindent (1) $X^{\prime}$ has the structure of a Mori fibre space $X^{\prime}\to Y$ , $K_{X^{\prime}}+B^{\prime}+A^{\prime}$ is the pullback of a $\Q$-divisor from $Y$, and $\mathrm{Supp}\ A^{\prime}$ dominates $Y$, or

\noindent (2) $K_{X^{\prime}}+B^{\prime}+(1-\varepsilon)A^{\prime}$ is nef for any sufficiently small rational $\varepsilon>0$.

\noindent Moreover, $K_{X^{\prime}}+B^{\prime}+A^{\prime}$ is semi-ample if and only if $K_{X}+B+A$ is semi-ample.
\elemma
\bproof We only need to prove that, if a $(K_{X}+B)$-MMP which is $(K_{X}+B+A)$-trivial terminates with a Mori fibre space $f:(X^{\prime},B^{\prime}+A^{\prime})\to Y$, then $\mathrm{Supp}\ A^{\prime}$ dominates $Y$. It is clear since $f$ only contracts curves which have positive intersections with $A^{\prime}$.
\eproof

We will use the following results on termination of flips.

\btheorem\label{Waldrontermination} (\cite[Theorem 1.6]{waldron2018lmmp} and \cite{hacon2022minimal}) Let $(X,B)$ be a projective lc threefold pair over a perfect field $k$ of characteristic $p>3$. If $M$ is an effective $\Q$-Cartier $\Q$-divisor on $X$, then any sequence of $(K_{X}+B)$-flips which are also $M$-flips terminates.
\etheorem

\blemma\label{terminaltermination}  Let $(X,B)$ be a $\Q$-factorial projective lc threefold pair over an algebraically closed field $k$ of characteristic $>3$ such that $K_{X}+B+A$ is nef.
If $X$ is terminal, then any $K_{X}$-MMP which is $(K_{X}+B)$-trivial terminates.
\elemma
\bproof Since every step of a $K_{X}$-MMP which is $(K_{X}+B)$-trivial is a step of a $K_{X}$-MMP, the assertion follows from \cite[Theorem 6.17]{kollar1998birational}.
\eproof

\subsection{Dlt modifications and adjunction}

The following result helps us to reduce some problems for lc pairs to $\Q$-factorial dlt pairs.

\btheorem\label{dltmodification}  Let $(X, B)$ be a lc threefold pair over a perfect field $k$ of characteristic $>3$.
Then $(X, B)$ has a crepant $\Q$-factorial dlt model. Moreover, we can modify $X$ so that it is terminal.
\etheorem
\bproof For the first assertion, see \cite[Theorem 1.6]{birkar2016existence} and \cite{hacon2022minimal}. Let us prove that we can make $X$ terminal. We take a crepant $\Q$-factorial dlt model $g: (X^{\prime},B^{\prime})\to (X,B)$ by the first assertion. Hence, by replacing $(X, B)$ by $(X^{\prime},B^{\prime})$, we may assume that $(X, B)$ is $\Q$-factorial and dlt. Let $U \subseteq X$ be the largest open set such that $(U, B|_{U})$ is a
snc pair. Then $\mathrm{codim}_{ X} (X\backslash U) \geq 2$. Let $f : (X^{\prime}, \Theta^{\prime})\to (X, 0)$ be a terminal model of $(X, 0)$ as in \cite[Theorem 1.7]{birkar2016existence} such that $K_{X^{\prime}} + \Theta^{\prime}= f^{\ast}K_{X}$. Then $f$ is an isomorphism over the smooth locus of $X$; in particular $f$ is an isomorphism
over $U$. Let $Z = X\backslash U$. Define 
$B^{\prime}:= \Theta^{\prime}+ f^{\ast}B$ on $X^{\prime}$
so that
$$K_{X^{\prime}} + B^{\prime} = f^{\ast}(K_{X}+ B),$$
and $(X^{\prime}, B^{\prime})$ is lc.

It remains to show that $(X^{\prime}, B^{\prime})$ is a dlt pair. Let $U^{\prime} = f^{-1}(U)$ and $Z^{\prime} = X^{\prime}\backslash U^{\prime}$.
Then $(U^{\prime}, B^{\prime}|_{U^{\prime}})$ is a snc pair. If $E$ is an exceptional divisor with center in $Z^{\prime}$,
then its center in $X$ is contained in $Z$. Hence $a(E, X^{\prime}, B^{\prime})= a(E, X, B) > -1$.
This completes the proof.
\eproof

For $\Q$-factorial dlt threefold pairs we have the following result on adjunction.

\btheorem\label{adjunction} Let $(X,B)$ be a $\Q$-factorial projective dlt threefold pair over a perfect field $k$ of characteristic $>0$. If $(K_{X}+ B)|_{\lfloor B\rfloor}$ is nef,  then $(K_{X}+ B)|_{\lfloor B\rfloor}$ is semi-ample.
\etheorem
\bproof By \cite[Remark 3.9]{hacon2022relative}, we know that all lc centres of $\Q$-factorial three-dimensional dlt pairs are normal up to a universal homeomorphism. Hence we can argue as in \cite[Section 5]{waldron2018lmmp} to prove that the $S_{2}$-fication (see \cite[2.2]{waldron2018lmmp} for example) of $\lfloor B\rfloor$ is a 
universal homeomorphism and $(K_{X}+ B)|_{\lfloor B\rfloor}$ is semi-ample.
\eproof

\subsection{Some known results on the abundance}

The following theorem collects the recent results towards the abundance conjecture in positive characteristics.

\btheorem\label{knownold} Let $(X,B)$ be a projective klt threefold pair over an algebraically closed field $k$ of characteristic $>5$ such that $K_{X}+B$ is nef. Assume that one of the following conditions holds:

\noindent (1) $\kappa(X,K_{X}+B)\geq 1$,
 
\noindent (2) the nef dimension $n(X,K_{X}+B)\leq 2$,

\noindent (3) the Albanese map $a_{X}: X\to\mathrm{Alb}(X)$ is non-trivial.
 
\noindent Then $K_{X}+ B$ is semi-ample. 
\etheorem
\bproof For (1) the case of $\kappa(X,K_{X}+B)=3$ is proved in \cite[Theorem 1.2]{birkar2017existence}, the case of $\kappa(X,K_{X}+B)=2$ is proved in \cite[Theorem 1.3]{waldron2017finite} and the case of $\kappa(X,K_{X}+B)=1$ is proved in \cite[Theorem 3.1]{zhang2019abundance} and \cite[Theorem A]{das2019abundance}. For (2) it is proved in \cite[Theorem 5]{witaszek2021canonical}. For (3) see \cite[Theorem 1.1]{zhang2020abundance} and \cite[Corollary 4.13]{witaszek2021canonical}. 
\eproof

Moreover, the non-vanishing theorem for terminal threefolds has been proved in \cite{xu2019nonvanishing}.

\btheorem (\cite[Theorem 1.1]{xu2019nonvanishing}) Let $X$ be a projective terminal threefold over an algebraically closed field $k$ of characteristic $>5$. If $K_{X}$ is pseudo-effective, then
$\kappa(X, K_{X})\geq 0$.
\etheorem

Based on it, the non-vanishing theorem for klt threefold pairs is proved in \cite{witaszek2021canonical}.

\btheorem\label{kltnonvanold} (\cite[Theorem 3]{witaszek2021canonical}) Let $(X,B)$ be a projective klt threefold pair over a perfect field $k$ of characteristic $>5$. If $K_{X}+B$ is pseudo-effective, then $\kappa(K_{X}+B)\geq 0$.
\etheorem

\section{Klt threefold pairs in characteristic $>3$}

In this section, we generalize the results in Subsection 2.6 to the case when the characteristic is greater than $3$. Note that in Subsection 2.6 we always assume that the characteristic of the base field is greater than $5$. Actually, the assumption of characteristic $>5$ is used for the following assertions.
Let $k$ be an algebraically closed field of characteristic $>5$. Then we have the following propositions hold:

\noindent \bf{P 1:} (MMP) We can run MMP for lc threefold pairs over $k$ (see \cite{hashizume2020minimal} for example).

\noindent \bf{P 2:} (Elliptic fibration) Let $g: X\to Z$ be a fibration of normal varieties of relative dimension one over $k$. Assume that the generic fiber $X_{\eta}$ of $g$ is a curve with arithmetic genus 
$p_{a}(X_{\eta}) = 1$. Then the geometric generic fiber $X_{\overline{\eta}}$ of $g$ is a
smooth elliptic curve over $\overline{K(Z)}$ (see \cite[Proposition 2.11]{zhang2020abundance}).

\noindent \bf{P 3:} (Dlt adjunction) Let $(X,B)$ be a $\Q$-factorial projective dlt threefold pair over $k$. Then every irreducible component of $\lfloor B \rfloor$ is normal.
If moreover $(K_{X}+ B)|_{\lfloor B\rfloor}$ is nef,  then it is semi-ample (see \cite[Section 2]{das2016adjunction} and \cite[Theorem 1.3]{waldron2018lmmp}).

\noindent \bf{P 4:} (Classification of surface $F$-singularity) Klt surface singularties over $k$ are strongly $F$-regular (see \cite{hara1998classification}).
\bremark These proposition are not independent. For example, the proof of \bf{P 1} uses \bf{P 4}.  
\eremark
Now we assume that the characteristic of $k$ is just greater than $3$. Then  \bf{P 1} and \bf{P 2} hold by \cite{hacon2022minimal} and \cite[Proposition 2.11]{zhang2020abundance}. Although \bf{P 3} may not hold, it is not far from being true.
More precisely, if $(X,B)$ is a $\Q$-factorial dlt threefold pair over $k$, then every irreducible component of $\lfloor B \rfloor$ is normal up to a universal homeomorphism by \cite[Remark 3.9]{hacon2022relative}. If, moreover, $(K_{X}+ B)|_{\lfloor B\rfloor}$ is nef,  then it is semi-ample by Theorem \ref{adjunction}. Finally, \bf{P 4} may not hold.

First, we generalize the results on subadditivity of Kodaira dimensions in \cite{zhang2020abundance} to the case when the characteristic is greater than $3$ (see Theorem \ref{subadditivity}). To do this, we need the following lemmas.

\blemma\label{replacelemma} (cf \cite[Lemma 4.10]{zhang2020abundance}) Let $(\hat{X},\hat{B})$ be a $\Q$-factorial  projective dlt threefold pair over an algebraically closed field $k$ of characteristic $>3$, and let $\hat{f}:\hat{X}\to Y$ be a fibration to a normal variety. Assume that $K_{\hat{X}}+\hat{B}$ is nef and 
$\hat{B}=G_{1}+G_{2}+\cdots +G_{n}$ is a sum of prime Weil divisors. Denote the normalization of  $G_{j}$ by $G^{\nu}_{j}$ for every $j=1,2,\cdots,n$. Then for every $j=1,2,\cdots,n$, $(K_{\hat{X}}+\hat{B})|_{G_{j}}$ is semi-ample.
Moreover,  a general fibre $F_{j}$ of the Iitaka fibration induced by $(K_{\hat{X}}+\hat{B})|_{G^{\nu}_{j}}$ is integral. We denote the image of $F_{j}$ along the normalization $G^{\nu}_{j}\to G_{j}$ by $\hat{F}_{j}$.

Assume in addition that

\noindent (a) there exist $N > 0$ and two different effective Cartier divisors $\hat{D}_{i}$
, $i= 1,2$ such that 
$$\hat{D}_{i}\sim N(K_{\hat{X}} + \hat{B})+ \hat{f}^{\ast}L_{i}$$
for some $L_{i}\in \mathrm{Pic}^{0}(Y)$ and that $\mathrm{Supp}\ \hat{D}_{i}\subseteq \mathrm{Supp}\ \hat{B}$,

\noindent (b) there exist effective divisors $\hat{G}_{1},\hat{G}_{2},\hat{G}^{\prime}_{1},\hat{G}^{\prime}_{2}$ such that 
$$\hat{D}_{1}=a_{11}\hat{G}_{1}+a_{12}\hat{G}_{2}+\hat{G}^{\prime}_{1},\hat{D}_{2}=a_{21}\hat{G}_{1}+a_{22}\hat{G}_{2}+\hat{G}^{\prime}_{2},$$ where $a_{11}> a_{21}\geq 0$ and $a_{22} > a_{12} \geq 0$, and

\noindent (c) $G_{1}, G_{2}$ are two irreducible components of $\hat{G}
_{1}, \hat{G}_{2}$ respectively, such that for $i, j\in \{1, 2\}$ and $i \neq j$, $\hat{F}_{j}$
dominates $Y$ and
$$\hat{F}_{j}\cap\mathrm{Supp}(\hat{G}^{\prime\prime}_{j}:=\hat{G}_{i}+\hat{G}^{\prime}_{1}+\hat{G}_{2}^{\prime})=\emptyset.$$

\noindent Then both $L_{1}$ and $L_{2}$ are torsion line bundles.

Furthermore, condition (c) holds, if for $j = 1, 2$, $G_{j}$ is not a component of $\hat{G}^{\prime\prime}_{j}$ and $\kappa(F_{j})\geq 0$.
\elemma
\bproof By Theorem \ref{adjunction}, we have $(K_{\hat{X}}+\hat{B})|_{\hat{B}}=(K_{\hat{X}}+\hat{B})|_{\lfloor\hat{B}\rfloor}$ is semi-ample. In particular, $(K_{\hat{X}}+\hat{B})|_{G_{j}}$, and hence $(K_{\hat{X}}+\hat{B})|_{G^{\nu}_{j}}$ are semi-ample for every $j=1,2,\cdots,n$.
Moreover, a general fibre $F_{j}$ of the Iitaka fibration induced by $(K_{\hat{X}}+\hat{B})|_{G^{\nu}_{j}}$ is integral by \cite[Proposition 2.1]{zhang2020abundance}. Hence, the first assertion holds.

Now we assume (a), (b) and (c). Note that 
$$(K_{\hat{X}}+\hat{B})|_{F_{1}}=\big( (K_{\hat{X}}+\hat{B})|_{G^{\nu}_{1}}\big)|_{F_{1}}\sim_{\Q}0$$ 
since $(K_{\hat{X}}+\hat{B})|_{G^{\nu}_{1}}$ is semi-ample and  $F_{1}$ is a general fibre of the Iitaka fibration of $(K_{\hat{X}}+\hat{B})|_{G^{\nu}_{1}}$. 
We have
\beq\begin{split}
a_{21}\hat{f}^{\ast}L_{1}|_{F_{1}}&\sim_{\Q} a_{21}(N(K_{\hat{X}} + \hat{B})+ \hat{f}^{\ast}L_{1})|_{F_{1}} \\
&\sim_{\Q} a_{21}\hat{D}_{1}|_{F_{1}}\ \ \ \ \ \ \ \ \ \ \ \ \ \ \ \ \ \ \ \ \ \ \ \ \ \ \ \ \ \ \ \ \ \ \ \ \mathrm{(by\ (a))} \\
&\sim_{\Q}a_{21}(a_{11}\hat{G}_{1}+a_{12}\hat{G}_{2}+\hat{G}^{\prime}_{1})|_{F_{1}}\ \ \ \mathrm{(by\ (b))}\\
&\sim_{\Q}a_{11}a_{21}\hat{G}_{1}|_{F_{1}}\ \ \ \ \ \ \ \ \ \ \ \ \ \ \ \ \ \ \ \ \ \ \ \ \ \ \ \ \ \ \ \mathrm{(by\ (c)).}   \notag
\end{split}\eeq
Similarly, $a_{11}\hat{f}^{\ast}L_{2}|_{F_{1}}\sim_{\Q}a_{11}a_{21}\hat{G}_{1}|_{F_{1}}$. Hence, we have $a_{21}\hat{f}^{\ast}L_{1}|_{F_{1}}\sim_{\Q}a_{11}\hat{f}^{\ast}L_{2}|_{F_{1}}$.
It follows that $a_{21}L_{1}\sim_{\Q}a_{11}L_{2}$ by \cite[Lemma 2.4]{zhang2020abundance}. Similarly, we have $a_{22}L_{1}\sim_{\Q}a_{12}L_{2}$. We then deduce that $L_{1}\sim_{\Q}L_{2}\sim_{\Q}0$ since $a_{11}> a_{21}\geq 0$ and $a_{22} > a_{12} \geq 0$. Hence the second assertion holds.

It remains to prove the third assertion. As $\kappa(F_{j})\geq 0$, we have the canonical
divisor $K_{F^{\nu}_{j}}\geq 0$, where $F_{j}^{\nu}$ is the normalization of $F_{j}$. Applying the adjunction formula, we get
\beq\begin{split}
0\sim_{\Q} (K_{\hat{X}} + \hat{B})|_{F_{j}^{\nu}}&\sim_{\Q}((K_{\hat{X}} + \hat{B})|_{G^{\nu}_{j}})|_{F_{j}^{\nu}} \\
&\sim_{\Q}((K_{\hat{X}}+G_{j})|_{G^{\nu}_{j}} + (\hat{B}-G_{j})|_{G^{\nu}_{j}})|_{F_{j}^{\nu}} \\
&\sim_{\Q} (K_{G^{\nu}_{j}}+C_{j})|_{F^{\nu}_{j}} + (\hat{B}-G_{j})|_{F_{j}^{\nu}}\\
&\sim_{\Q}K_{F^{\nu}_{j}}+C_{j}|_{F^{\nu}_{j}} + (\hat{B}-G_{j})|_{F_{j}^{\nu}}\notag
\end{split}\eeq
where $C_{j}\geq 0$ on $G^{\nu}_{j}$. It implies that $(\hat{B}-G_{j})|_{F_{j}^{\nu}}\leq 0$.
Since $F_{j}$ is general,  $\hat{F}_{j}$ is not contained in $\hat{B}-G_{j}$. Hence, 
$\hat{F}_{j}\cap\mathrm{Supp}(\hat{B}-G_{j})=\emptyset$.
By our assumption, $G_{j}$ is not a component of $\hat{G}^{\prime\prime}_{j}$. Thus, 
$\mathrm{Supp}(\hat{G}^{\prime\prime}_{j})\subseteq\mathrm{Supp}(\hat{B}-G_{j})$.
It follows that $\hat{F}_{j}\cap\mathrm{Supp}(\hat{G}^{\prime\prime}_{j})=\emptyset$.
\eproof

\blemma\label{kodairadimsemi-ample} Let $(X,B)$ be a projective klt threefold pair over an algebraically closed field $k$ of characteristic $>3$. Assume that $K_{X}+B$ is nef and  $\kappa(X,K_{X}+B)\geq 1$. Then $K_{X}+B$ is semi-ample.
\bproof The case of $\kappa(X,K_{X}+B)=3$ follows from \cite[Theorem 1.3]{hacon2022minimal}.   In the cases of $\kappa(X,K_{X}+B)= 1\ \mathrm{or}\ 2$, 
the assertion is proved when the characteristic of $k$ is greater than $5$ in \cite[Theorem 1.3]{waldron2017finite}, \cite[Theorem 3.1]{zhang2019abundance} and \cite[Theorem A]{das2019abundance}. And it uses the assumption of characteristic $>5$ for \bf{P 1}.  When the characteristic of $k$ is greater than $3$, by Theorem \ref{lcmmp}, \bf{P 1} also holds. Hence we 
 can argue as in the proofs of \cite[Theorem 1.3]{waldron2017finite}, \cite[Theorem 3.1]{zhang2019abundance} and \cite[Theorem A]{das2019abundance} to prove the assertion.
\eproof
\elemma

Now we can deduce the following result on subadditivity of Kodaira dimensions in characteristic $>3$.

\btheorem\label{subadditivity} Let $f: X\to Y$ be a fibration from a $\Q$-factorial
projective threefold to a smooth projective variety of dimension $1$ or $2$, over an algebraically closed field $k$ of characteristic $p > 3$. Assume that there is an effective $\Q$-divisor $B$ on
$X$ such that $(X, B)$ is klt. Assume that $Y$ is of maximal Albanese dimension. Moreover, we assume that if $\kappa(X_{\eta}, K_{X_{\eta}}+ B_{\eta}) =\mathrm{dim}\ X-\mathrm{dim}\ Y- 1,$ where $X_{\eta}$ is the generic fibre of $f$ and $K_{X_{\eta}}+ B_{\eta}:=(K_{X}+B)|_{X_{\eta}}$, then $B$ does not intersect the generic fibre $X_{\xi}$ of the relative Iitaka fibration $I: X \ratmap Z$ induced by $K_{X}+ B$ on $X$ over $Y$.

\noindent Then 
$$\kappa(X,K_{X}+B)\geq \kappa(X_{\eta},K_{X_{\eta}}+ B_{\eta})+\kappa(Y).$$
\etheorem
\bproof The case when the characteristic is greater than $5$ is proved in \cite[Theorem 1.4]{zhang2020abundance}. Using Theorem \ref{lcmmp} and Lemma \ref{kodairadimsemi-ample} we can argue as in the proof of \cite[Theorem 1.4]{zhang2020abundance} except in the cases when

\noindent (1) $Y$ is an elliptic curve or a simple abelian surface, and $K_{X}+B$ is $f$-big, or

\noindent (2) $Y$ is an elliptic curve, $\kappa(X_{\eta}, K_{X_{\eta}}+ B_{\eta}) =1$ and $B$ does not intersect the generic fibre $X_{\xi}$ of the relative Iitaka fibration $I: X \ratmap Z$ induced by $K_{X}+ B$ on $X$ over $Y$.

Now we assume that we are in one of these cases. 
We first make some reductions as follows. In the case (1), if the characteristic of $k$ is greater than $5$, then the proof of \cite[Theorem 4.2]{zhang2020abundance} reduces the assertion to the case when 

\noindent $\bullet$ the denominators of coefficients of $B$ are not divisible by $p$,

\noindent $\bullet$ $K_{X}+B$ is a nef and $f$-ample,

\noindent $\bullet$ $\nu(K_{X}+B)\leq 2$,

\noindent $\bullet$ there exist a sufficiently divisible positive integer $l$ and a coherent sheaf $\sF$ such that $\sF$ is a subsheaf of $f_{\ast}\sO_{X}(l(K_{X}+B))$,

\noindent $\bullet$ there exists an isogeny $\tau:Y_{1}\to Y$ between abelian varieties, some $P_{i}\in \mathrm{Pic}^{0}(Y_{1})$ and a generically surjective homomorphism 
$$\tau^{\ast}\sF\cong \oplus_{i=1}^{r_{1}}P_{i}.$$

\noindent In the case (2), if the characteristic of $k$ is greater than $5$, the proof of \cite[Theorem 4.3]{zhang2020abundance} reduces the assertion to the case when 

\noindent $\bullet$ $K_{X}+B$ is nef,

\noindent $\bullet$ there exists a commutative diagram
\beq\begin{tikzcd}
W\arrow[r,"\sigma "]\arrow[d,swap,"h"]&X\arrow[d,"f"] \\
Z\arrow[r,"g"] &Y\notag
\end{tikzcd}
\eeq
where $\sigma$ is a log resolution, $h$ is a fibration to a smooth projective surface which is birational to the relative Iitaka fibration induced by $\sigma^{\ast}(K_{X}+B)$ on $W$ over $Y$,

\noindent $\bullet$ there exists a nef and $g$-big divisor $C$ on $Z$ such that $\sigma^{\ast}(K_{X}+B)\sim_{\Q}h^{\ast}C$,

\noindent $\bullet$ the geometric generic fibre of $g$ is either a smooth elliptic curve or a rational curve,

\noindent $\bullet$ $\nu(Z,C)=1$,

\noindent $\bullet$ there exist a sufficiently divisible positive integer $l$ and a nef sub-vector bundle $V$ of $f_{\ast}\sO_{X}(l(K_{X}+B))$ of rank $r\geq 2$,

\noindent $\bullet$ there exists a flat base change $\pi:Y_{2}\to Y$ between elliptic curves such that 
$$\pi^{\ast}V\cong \oplus_{i=1}^{r_{2}}L_{i}^{\prime},$$ where $L_{i}^{\prime}\in \mathrm{Pic}^{0}(Y_{2}).$

\noindent When the characteristic of $k$ is greater than $3$, using Theorem \ref{lcmmp} and Lemma \ref{kodairadimsemi-ample} we can also argue as in the proofs of \cite[Theorem 4.2 and Theorem 4.3]{zhang2020abundance} to make such reductions.

If the characteristic of $k$ is greater than $5$, then the argument in \cite[Step 2,3 of the proof of Theorem 4.2 and Step 2,3 of the proof of Theorem 4.3]{zhang2020abundance} implies that there exist an integer $m_{1}$ and some divisors $D_{i}\in |m_{1}(K_{X}+B)+f^{\ast}L_{i}|$, $i=1,2,\cdots,r$ for some $L_{i}\in \mathrm{Pic}^{0}(Y)$. Moreover, we can construct a pair $(\hat{X},\hat{B})$ and divisors $\hat{D}_{1},\hat{D}_{2}$ satisfying all conditions of Lemma \ref{replacelemma}. When the characteristic of $k$ is greater than $3$, using Theorem \ref{lcmmp}, we can also argue as in the proofs of \cite[Theorem 4.2 and Theorem 4.3]{zhang2020abundance} to prove these assertions.
By Lemma \ref{replacelemma}, $L_{1}$ and $L_{2}$ are torsions. Hence there exist a sufficiently divisible integer $N>0$ and two different divisors among $D_{i}$, say, $D_{1}\neq D_{2}$ such that
$$ND_{j}\in |Nm_{1}(K_{X}+B)+NL_{j}|=|Nm_{1}(K_{X}+B)|$$
for $j=1,2$. Hence we have $\kappa(X,K_{X}+B)\geq 1$. In the case (2), it implies that
$$\kappa(X,K_{X}+B)\geq 1=\kappa(X_{\eta},K_{X_{\eta}}+B_{\eta}).$$ In the case (1), by Lemma \ref{kodairadimsemi-ample}, $K_{X}+B$ is semi-ample. Thus for a sufficiently divisible $M>0$, the linear system $|M(K_{X}+B)|$ has no base point. Since $K_{X_{\eta}}+B_{\eta}$ is big, the restriction $|M(K_{X}+B)||_{X_{\eta}}$ on the generic fibre $X_{\eta}$ defines a generically finite morphism. It implies that
$$\kappa(X,K_{X}+B)\geq \mathrm{dim}\ X_{\eta}= \kappa(X_{\eta},K_{X_{\eta}}+B_{\eta}).$$ 
In conclusion, the assertion holds.
\eproof

Using this result on subadditivity of Kodaira dimensions in characteristic $>3$, we deduce the  following results on the abundance with non-trivial Albanese maps in characteristic $>3$.

\blemma\label{nonunirulednontrivialAlb} Let $(X,B)$ be a $\Q$-factorial projective klt threefold pair over an algebraically closed field $k$ of characteristic $> 3$. Assume that $K_{X}+B$ is nef, $X$ is non-uniruled and the Albanese map $a_{X}: X\to \mathrm{Alb}(X)$ is non-trivial. Then $K_{X}+B$ is semi-ample.
\elemma
\bproof The case when the characteristic of $k$ is greater than $5$ is proved in \cite[Theorem 1.1]{zhang2019abundance}. When the characteristic of $k$ is greater than $3$, by the proof of \cite[Theorem 1.1]{zhang2019abundance}, we only need to prove the following assertions.

\noindent (1) Let $f_{1}: X_{1}\to Y_{1}$ be a separable fibration from a
smooth projective threefold to a smooth projective variety of dimension $1$ or $2$ over $k$. Denote by $\tilde{X}_{1,\overline{\eta}}$ a smooth projective birational model of $X_{1,\overline{\eta}}$, where $X_{1,\overline{\eta}}$ is the geometric generic fibre of $f_{1}$. Then
$$\kappa(X_{1})\geq \kappa(\tilde{X}_{1,\overline{\eta}})+\kappa(Y_{1}).$$

\noindent (2)  Let $X_{2}$ be a $\Q$-factorial projective klt threefold over $k$ with $K_{X_{2}} \sim_{\Q} 0$, and let $D$ be an effective and nef $\Q$-divisor on $X_{2}$. Assume that $X_{2}$ has a morphism $f_{2}: X_{2}\to Y_{2}$ to an elliptic curve and that $X_{2,\overline{\eta}}$ has at most canonical singularities, where $X_{2,\overline{\eta}}$ is the geometric generic fibre of $f_{2}$. Then either $D = 0$ or $\kappa(X_{2}, D)\geq 1$.

(1) is proved when the characteristic of $k$ is greater than $5$ in \cite[Corollary 2.9]{zhang2019abundance}. It uses the assumption of characteristic $>5$ for the fact that canonical singularities over $k$ are $F$-pure. This fact holds in characteristic $5$ by \cite[Theorem 1.2]{hara1998classification}. Hence (1) follows from the proof of \cite[Corollary 2.9]{zhang2019abundance}. For (2), it suffices to show that if $\kappa(X_{2}, D)=0$, then $D=0$. We assume that $\kappa(X_{2}, D)=0$. We denote the generic fibre of $f_{2}$ by $X_{2,\eta}$.
Note that 
$$D_{\eta}:=D|_{X_{2,\eta}}\sim_{\Q}K_{X_{2,\eta}}+ D_{\eta}$$ 
and $(X_{2,\eta},D_{\eta})$ is lc after replacing $D$ by a small multiple.
By Theorem \ref{surfaceabundance}, $D_{\eta}$ is semi-ample. Hence $\kappa(X_{2,\eta},D_{\eta})\geq 0$. If $\kappa(X_{2,\eta},D_{\eta})\neq 1$, then by Theorem \ref{subadditivity}, we have $\kappa(X_{2,\eta},D_{\eta})=0$. Hence $D_{\eta}\sim_{\Q}0$. Note that $f_{2}$ is equidimensional since $Y_{2}$ is a normal curve. By Lemma \ref{descend}, $D$ descends to an effective $\Q$-divisor on $Y_{2}$. Hence $D=0$. Otherwise, we have $\kappa(X_{2,\eta},D_{\eta})= 1$. Then we may apply the proof of \cite[Corollary 2.10]{zhang2019abundance} to the case of the characteristic of $k$ is greater than $3$. Therefore, the assertion holds.
\eproof
\bremark The non-uniruled assumption is used in the proof of \cite[Theorem 1.1]{zhang2019abundance}.
\eremark

\btheorem\label{nontrivialAlb} Let $(X,B)$ be a $\Q$-factorial projective klt threefold pair over an algebraically closed field $k$ of characteristic $> 3$. Assume that $K_{X}+B$ is nef and the Albanese map $a_{X}: X\to \mathrm{Alb}(X)$ is non-trivial. Denote by $f: X\to Y$ the fibration arising from the Stein factorization of $a_{X}$ and by $X_{\eta}$ the generic fiber of $f$. Assume moreover
that $B = 0$ if

\noindent (1) $\mathrm{dim}\ Y=2$ and $\kappa(X_{\eta},(K_{X}+B)|_{X_{\eta}})=0$, or

\noindent (2) $\mathrm{dim}\ Y=1$ and $\kappa(X_{\eta},(K_{X}+B)|_{X_{\eta}})=1$.

\noindent Then $K_{X}+B$ is semi-ample.
\etheorem
\bproof The case when the characteristic is greater than $5$ is proved in \cite[Theorem 1.2]{zhang2020abundance}.
By Lemma \ref{nonunirulednontrivialAlb} we can assume that $X$ is uniruled.  Moreover, by Lemma \ref{kodairadimsemi-ample}, we can assume that $\kappa(X,K_{X}+B)\leq 0$.

Since $X$ is uniruled, we have $\mathrm{dim}\ Y=1\ \mathrm{or}\ 2$. Note that $K_{X_{\eta}}+B_{\eta}$ is semi-ample by the abundance for surfaces (Theorem \ref{surfaceabundance}) and curves. In particular, $\kappa(X_{\eta},K_{X_{\eta}}+B_{\eta})\geq 0$. Therefore by Theorem \ref{subadditivity}, we have $\kappa(X,K_{X}+B)=0$, and hence $\kappa(Y)=\kappa(X_{\eta},K_{X_{\eta}}+B_{\eta})=0$.
If $\mathrm{dim}\ Y=1$, then the assertion is proved when the characteristic of $k$ is greater than $5$ in \cite[Theorem 4.4]{zhang2020abundance}.
Using Theorem \ref{subadditivity} we can argue as in the proof of
\cite[Theorem 4.4]{zhang2020abundance} to prove that $K_{X}+B$ is semi-ample. Otherwise, we have $\mathrm{dim}\ Y=2$.
Then $B=0$ by our assumption and $f$ is an elliptic fibration by \cite[Proposition 2.11]{zhang2020abundance}. Hence $X$ is non-uniruled. We obtain a contradiction.
Thus, $K_{X}+B$ is semi-ample.
\eproof
\bcorollary\label{terminalnonvan} Let $X$ be a projective terminal threefold over an algebraically closed field $k$ of characteristic $>3$. If $K_{X}$ is pseudo-effective, then
$\kappa(X, K_{X})\geq 0$.
\ecorollary
\bproof The case when the characteristic of $k$ is greater than $5$ is proved in \cite[Theorem 1.1]{xu2019nonvanishing}. Using Theorem \ref{nontrivialAlb}, we can argue as in the proof of \cite[Theorem 1.1]{xu2019nonvanishing} to prove the assertion.
\eproof

Now we can generalize Theorem \ref{knownold} to the case when the characteristic is greater than $3$.

\btheorem\label{known} Let $(X,B)$ be a projective klt threefold pair over an algebraically closed field $k$ of characteristic $>3$ such that $K_{X}+B$ is nef. Assume that one of the following conditions holds:

\noindent (1) $\kappa(X,K_{X}+B)\geq 1$,
 
\noindent (2) the nef dimension $n(X,K_{X}+B)\leq 2$,

\noindent (3) the Albanese map $a_{X}: X\to\mathrm{Alb}(X)$ is non-trivial.
 
\noindent Then $K_{X}+ B$ is semi-ample. 
\etheorem
\bproof See Lemma \ref{kodairadimsemi-ample} for (1). For (2),
it is proved when the characteristic of $k$ is greater than $5$ in \cite[Theorem 5]{witaszek2021canonical}.
Using Theorem \ref{nontrivialAlb} in the case of $n(X,K_{X}+B)=0$, we can argue as in the proof of \cite[Theorem 5]{witaszek2021canonical} to prove the assertion.
For (3), it is proved when the characteristic of $k$ is greater than $5$ in \cite[Corollary 4.13]{witaszek2021canonical}.
Using Theorem \ref{nontrivialAlb} and (2), we can argue as in the proof of \cite[Corollary 4.13]{witaszek2021canonical} to prove the assertion.
\eproof

Moreover, we can deduce the non-vanishing theorem for klt threefold pairs in characteristic $>3$.

\btheorem\label{kltnonvan} Let $(X,B)$ be a projective klt threefold pair over an algebraically closed field $k$ of characteristic $>3$. If $K_{X}+B$ is pseudo-effective, then $\kappa(K_{X}+B)\geq 0$.
\etheorem
\bproof It is proved when the characteristic of $k$ is greater than $5$ in \cite[Theorem 3]{witaszek2021canonical}.
Using Corollary \ref{terminalnonvan} and (2) of Theorem \ref{known}, we can argue as in the proof of \cite[Theorem 3]{witaszek2021canonical} to prove the assertion.
\eproof

\section{Nonvanishing theorem for lc threefold pairs}

In this section we show the nonvanishing theorem for projective lc threefold pairs.
First, we recall a standard lemma on modifying a pair by some birational transform.

\blemma\label{standardlemma} Let $(X,B)$ be a $\Q$-factorial dlt threefold pair over an algebraically closed field $k$ of characteristic $>3$. Suppose that $K_{X}+B$ is nef and there exists an effective $\Q$-divisor $D$ such that
$D\equiv K_{X}+B$. Then there exists a $\Q$-factorial dlt pair $(Y,B_{Y})$ such that

\noindent (1) $K_{Y}+B_{Y}$ is nef,

\noindent (2) $n(K_{Y}+B_{Y})=n(K_{X}+B)$,

\noindent (3) $\kappa(K_{X}+B)\leq \kappa(K_{Y}+B_{Y})\leq \kappa(K_{X}+B+rD)$ for some $r>0$,

\noindent (4) $K_{Y}+B_{Y}\equiv \Delta$ for an effective $\Q$-divisor $\Delta$ with $\mathrm{Supp}\ \Delta\subseteq \lfloor B_{Y}\rfloor$,

\noindent (5) $(Y\backslash\mathrm{Supp}\ \Delta,B_{Y})\cong (X\backslash\mathrm{Supp}\  D,B)$.

\noindent Moreover, if $D \sim_{\Q} K_{X}+B$, then $K_{Y}+ B_{Y} \sim_{\Q} \Delta$ in (4).
\elemma
\bproof It follows from Theorem \ref{lcmmp} and the proof of \cite[Lemma 4.6]{witaszek2021canonical}.
\eproof

The following lemma is proved by Witaszek via his weak canonical bundle formula.

\blemma\label{specialweakcan} (\cite[Lemma 4.8]{witaszek2021canonical}) Let $(X,B)$ be a projective $\Q$-factorial threefold pair over an algebraically closed field $k$ of characteristic $>3$ such that the coefficients of $B$ are at most one. Assume that $L := K_{X}+B$ is nef and $n(L) = 2$. Then the following hold:

\noindent (1) there exists an effective $\Q$-divisor $D$ such that $L \equiv D$,

\noindent (2) if $L|_{\mathrm{Supp}\ D} \sim_{\Q} 0$ for some $D$ as above, then $\kappa(L) \geq 0$,

\noindent (3) if $L|_{\mathrm{Supp}\ D}\not\equiv 0$ for some $D$ as above, or $L|_{\mathrm{Supp}\ D} \sim_{\Q}0$ and $L\sim_{\Q} D$, then $\kappa(L)= 2$.
\elemma

Then we can deduce the following proposition.

\bproposition\label{weakcan} Let $(X,B)$ be a projective lc threefold pair over an algebraically closed field $k$ of characteristic $>3$. If $K_{X}+B$ is nef and $n(X,K_{X}+B)= 2$, then $\kappa(K_{X}+B)=2$ .
\eproposition
\bproof The proof is similar to the proof of \cite[Proposition 4.10]{witaszek2021canonical}.
By Theorem \ref{dltmodification} replacing $(X,B)$ by a $\Q$-factorial dlt model,  we may assume that $(X,B)$ is $\Q$-factorial and dlt. By Lemma \ref{specialweakcan}, there exists an effective $\Q$-divisor $D$ satisfying $K_{X}+B\equiv D$. Now by Lemma \ref{standardlemma} we have a $\Q$-factorial dlt pair
$(Y,B_{Y})$ such that for some $r > 0$,

\noindent $\bullet$ $K_{Y}+B_{Y}$ is nef,

\noindent $\bullet$ $n(K_{Y}+B_{Y}) = n(K_{X}+B)$ and $\kappa(K_{Y}+B_{Y})\leq \kappa(K_{X} +B+rD)$,

\noindent $\bullet$ $K_{Y}+B_{Y} \equiv E_{Y}$, where $E_{Y}$ is an effective $\Q$-divisor such that $\mathrm{Supp}\ E_{Y}\subseteq \lfloor B_{Y}\rfloor$.

\noindent By Theorem \ref{adjunction}, $(K_{Y}+B_{Y})|_{\lfloor B_{Y}\rfloor}$ , and hence $(K_{Y}+B_{Y})|_{\mathrm{Supp}\ E_{Y}}$ are semi-ample. Applying Lemma \ref{specialweakcan}  to $(Y ,B_{Y})$ and $E_{Y}$ , we have $\kappa(K_{Y} +B_{Y})\geq 0$.

We claim that in fact $\kappa(K_{Y}+B_{Y})\geq 2$. We apply Lemma \ref{standardlemma} to $(Y,B_{Y})$
and an effective $\Q$-divisor which is $\Q$-linearly equivalent to $K_{Y}+ B_{Y}$ ,  then we obtain a $\Q$-factorial dlt pair $(Z,B_{Z})$ satisfying

\noindent $\bullet$ $K_{Z}+B_{Z}$ is nef,

\noindent $\bullet$ $n(K_{Z}+B_{Z}) = n(K_{Y}+B_{Y})$ and $\kappa(K_{Z}+B_{Z})= \kappa(K_{Y} +B_{Y})$,

\noindent $\bullet$ $K_{Z}+B_{Z} \sim_{\Q} E_{Z}$, where $E_{Z}$ is an effective $\Q$-divisor such that $\mathrm{Supp}\ E_{Y}\subseteq \lfloor B_{Y}\rfloor$.

Similarly, we have $(K_{Z}+B_{Z})|_{\mathrm{Supp}\ E_{Z}}$ is semi-ample. Therefore, by Lemma \ref{specialweakcan} we have $\kappa(K_{Y}+B_{Y})=\kappa(K_{Z}+B_{Z})=2$. It implies that $\kappa(K_{X}+B+rD)\geq 2$. Since $K_{X}+B\equiv D$, it is clear that $(K_{X}+B)|_{D}\not\equiv 0$. Finally, by Lemma \ref{specialweakcan} we have $\kappa(K_{X}+B)=2$.
\eproof

Now we can prove the nonvanishing theorem for projective lc threefold pairs.

\btheorem\label{nonvan} Let $(X,B)$ be a projective lc threefold pair over a perfect field $k$ of characteristic $>3$. If $K_{X}+ B$ is pseudo-effective, then $\kappa(X,K_{X}+B)\geq 0$.
\etheorem
\bproof We pass to an uncountable algebraically closed field. Replacing $(X,B)$ by its log minimal model by Theorem \ref{lcmmp} , we can assume that $K_{X}+B$ is nef. 
By Theorem \ref{dltmodification}, we can take a $\Q$-factorial dlt model $(X^{\prime},B^{\prime})$ of $(X,B)$ such that $(X^{\prime},B^{\prime})$ is $\Q$-factorial and dlt, and moreover $X^{\prime}$ is terminal. We replace $(X,B)$ by $(X^{\prime},B^{\prime})$. If $\lfloor B\rfloor =0$, then the proposition follows from Theorem \ref{kltnonvan}. Hence we can assume that $\lfloor B\rfloor \neq 0$.

Now by Definition \ref{MMPKtrivial} we run a $K_{X}$-MMP which is $(K_{X}+B)$-trivial. By Lemma \ref{terminaltermination}, it terminates with a pair $(X^{\prime\prime},B^{\prime\prime})$. Note that $(X,(1-\varepsilon)B)$ is klt and every step of a $K_{X}$-MMP which is $(K_{X}+B)$-trivial is a step of a $(K_{X}+(1-\varepsilon)B)$-MMP for any sufficiently small rational $\varepsilon>0$. Hence we have $(X^{\prime\prime},(1-\varepsilon)B^{\prime\prime})$ is klt for any sufficiently small rational $\varepsilon>0$.
If $K_{X^{\prime\prime}}+(1-\varepsilon)B^{\prime\prime}$ is nef for any sufficiently small rational $\varepsilon>0$, then we have $\kappa(K_{X^{\prime\prime}}+(1-\varepsilon)B^{\prime\prime})\geq 0$ by Theorem \ref{kltnonvan}
since $(X^{\prime\prime},(1-\varepsilon)B^{\prime\prime})$ is klt.
Hence we have 
$$\kappa(K_{X}+B)=\kappa(K_{X^{\prime\prime}}+B^{\prime\prime})\geq\kappa(K_{X^{\prime\prime}}+(1-\varepsilon)B^{\prime\prime})\geq 0.$$ 

Otherwise, by Lemma \ref{Output} we get a Mori fibre space
\beq\begin{tikzcd}
X\arrow[r,dashrightarrow]& X^{\prime\prime}\arrow[d,"f"]\\
&Z\notag
\end{tikzcd}\eeq
and $\Q$-divisors $C$ on $Z$ such that
$$K_{X^{\prime\prime}}+B^{\prime\prime}\sim_{\Q} f^{\ast}C.$$
Hence we have 
$$n(K_{X}+B)\leq \mathrm{dim}\ Z \leq 2.$$ 
If $n(K_{X}+B)= 2$, by Proposition \ref{weakcan} we have $\kappa(K_{X}+B)=2$. 
If $n(K_{X}+B)= 1$, then by Theorem \ref{nefreductionmap} we get a nef reduction map of $K_{X}+B$, $g: X\to Z^{\prime}$. Then $g$ is an equidimensional morphism since $Z^{\prime}$ is a normal curve and $g$ is proper over the generic point of $Z^{\prime}$.  By Theorem \ref{surfaceabundance} we have $(K_{X}+B)|_{G}\sim_{\Q}0$, where $G$ is the generic fibre of $g$. Hence by Lemma \ref{descend}, $K_{X}+B$ descends to an ample divisor on $Z^{\prime}$. Therefore $K_{X}+B$ is semi-ample.

If $n(K_{X}+B)= 0$, then $K_{X}+B$ is numerically trivial.
By Theorem \ref{lcmmp}, there exists a $(K_{X}+B-\lfloor B\rfloor)$-MMP which terminates. Since $\lfloor B\rfloor > 0$, this MMP terminates with a Mori fibre space
\beq\begin{tikzcd}
X\arrow[r,dashrightarrow]& Y\arrow[d,"f^{\prime}"]\\
&Z^{\prime\prime}\notag
\end{tikzcd}\eeq
There are $\Q$-divisors $C^{\prime}$ on $Z^{\prime\prime}$, $B_{Y}$ on $Y$ such that $B_{Y}$ is the birational transform of $B$ on $Y$ and
$$K_{Y}+B_{Y}\sim_{\Q} f^{\prime\ast}C^{\prime}.$$
Now by Theorem \ref{dltmodification} we can take a dlt modification 
$$\mu:(Y^{\prime},B_{Y^{\prime}})\to (Y,B_{Y}).$$
Note that $\lfloor B_{Y^{\prime}}\rfloor$ dominates $Z^{\prime\prime}$  since $f^{\prime}$ only contract curves which have positive intersections with $\lfloor B_{Y}\rfloor$. Since $(K_{Y^{\prime}}+B_{Y^{\prime}})|_{\lfloor B_{Y^{\prime}} \rfloor}$ is semi-ample by Theorem \ref{adjunction}, we deduce that $C^{\prime}$, and hence $K_{X}+B$ are semi-ample by Lemma \ref{pullbacksemi-ample}.
\eproof

As a corollary, we have the following result on termination of flips.

\btheorem\label{termination} Let $(X,B)$ be a projective lc threefold pair defined over a perfect field $k$ of characteristic $p > 3$ such that $K_{X}+B$ is pseudo-effective. Then every sequence of $(K_{X}+B)$-flips terminates. In particular, any $(K_{X}+B)$-MMP terminates with a minimal model.
\etheorem
\bproof By Theorem \ref{nonvan}, we have $\kappa(K_{X}+B)\geq 0$. Then the proposition follows from Theorem \ref{Waldrontermination}. 
\eproof

\section{Abundance conjecture for lc threefold pairs}
In this section we show the abundance for lc threefold pairs whose Kodaira dimension $\geq 1$. To be precise, we prove the following result.

\btheorem\label{mainthm} Let $(X,B)$ be a projective lc threefold pair over an algebraically closed field  $k$ of characteristic $>3$. If $K_{X}+B$ is nef and $\kappa(X,K_{X}+B)\geq 1$, then $K_{X}+B$ is semi-ample.
\etheorem

\subsection{Preparation}

Before proving Theorem \ref{mainthm}, we make some preparations.

\blemma\label{kodnef} Let $X$ be a normal projective variety of dimension $3$ over an algebraically closed field, and $D$ be a nef $\Q$-Cartier $\Q$-divisor on $X$ such that $\kappa(X, D) = 2$.
Then $n(X,D)= 2$.
\elemma
\bproof We pass to an uncountable algebraically closed field. Consider the Iitaka map of $D$. After resolving the indeterminacies and replacing $D$ by its pullback, we can assume that the Iitaka map of $D$ is a morphism. Since $D$ is nef and not big, it has to be numerically trivial on all fibres of the Iitaka map. Hence we have $n(X,D)\leq 2$. Then by the equality $\kappa(X,D)\leq n(X,D)$ we have $n(X,D)=2$.
\eproof

\blemma\label{diagram} Let $X$ be a normal projective variety of dimension $3$ over an
uncountable algebraically closed field of characteristic $>0$. Assume $D$ is a nef $\Q$-Cartier $\Q$-divisor on $X$ such that $\kappa(X, D) = 2$.
Then $D$ is endowed with a map $h: X\to Z$ to a normal proper algebraic space $Z$ of dimension $2$.

If moreover $D|_{G}\sim_{\Q}0$, where $G$ is the generic fibre of $h$,
then there exists a commutative diagram
\beq\begin{tikzcd}
X_{1}\arrow[r,"\phi "]\arrow[d,swap,"h_{1}"]&X\arrow[d,"h"] \\
Z_{1}\arrow[r,"\psi"] &Z\notag
\end{tikzcd}
\eeq
where $Z_{1}$ is a smooth projective surface, $X_{1}$ is a normal projective threefold, $\phi,\psi$ are birational morphisms, and $h_{1}: X_{1}\to Z_{1}$ is an equidimensional fibration. Moreover, there exists a nef and big $\Q$-divisor $D_{1}$ on $Z_{1}$ such that $\phi^{\ast}D\sim_{\Q} h^{\ast}_{1}D_{1}$.
\elemma
\bproof By Lemma \ref{kodnef}, we have $\kappa(X, D)=n(X,D)= 2$. Hence by Lemma \ref{EWM}, $D$ is endowed with a map $h:X\to Z$ to a normal proper algebraic space $Z$ of dimension $2$. 

Assume moreover $D|_{G}\sim_{\Q}0$, where $G$ is the generic fibre of $h$. By Theorem \ref{nefreductionmap} we get a nef reduction map $f:X\dashrightarrow Y$ of $D$. Resolving the indeterminacies of $f$ and replacing $D$ by its pullback, we can assume that $f:X\to Y$ is a morphism to a normal surface.

Now we apply Lemma \ref{descend} to $f$ and $D$. Then we get a commutative diagram
\beq\begin{tikzcd}
X^{\prime}\arrow[d,swap,"f^{\prime} "]\arrow[r,"\phi^{\prime} "]&X\arrow[d,"f"] \\
Z^{\prime}\arrow[r,"\psi^{\prime} "] &Y\notag
\end{tikzcd}
\eeq
with $\phi^{\prime},\psi^{\prime}$ projective birational, and an $\Q$-divisor $C$ on $Z^{\prime}$ such that $\phi^{\prime\ast}D \sim_{\Q} f^{\prime\ast}C$. Moreover we can apply the flattening trick \cite[Theorem 5.2.2]{raynaud1971criteres} to $f^{\prime}$, and we get the following commutative diagram
\beq\begin{tikzcd}
X_{1}\arrow[d,swap,"h_{1}"]\arrow[r,"\phi^{\prime\prime} "]&X^{\prime}\arrow[d,swap,"f^{\prime} "]\arrow[r,"\phi^{\prime} "]&X\arrow[d,"f"] \\
Z_{1}\arrow[r,"\psi^{\prime\prime} "]&Z^{\prime}\arrow[r,"\psi^{\prime} "] &Y\notag
\end{tikzcd}
\eeq
where $Z_{1}$ is a normal projective surface, $X_{1}$ is a normal projective threefold, $\phi^{\prime\prime},\psi^{\prime\prime}$ are birational morphisms, and $h_{1}: X_{1}\to Z_{1}$ is a flat fibration.
Replacing $Z_{1}$ by a smooth resolution and $X_{1}$ by the normalization of main component of the fibre product of $h_{1}$ and the resolution, we may assume that $Z_{1}$ is smooth.

Let $\phi:=\phi^{\prime}\circ \phi^{\prime\prime}, D_{1}:=\psi^{\prime\prime\ast}C.$
Then we have 
$$\phi^{\ast}D\sim_{\Q} h^{\ast}_{1}D_{1}.$$
Since $h_{1}$ only contracts curves which are $\phi^{\ast}D$-numerically trivial, we know that the morphism $h\circ\phi: X_{1}\to Z$ factors through
$h_{1}$. In other words, there exists a natural map $\psi:Z_{1}\to Z$ making the following diagram commutative
\beq\begin{tikzcd}
X_{1}\arrow[r,"\phi "]\arrow[d,swap,"h_{1}"]&X\arrow[d,"h"] \\
Z_{1}\arrow[r,"\psi"] &Z\notag
\end{tikzcd}
\eeq
This completes the proof of the lemma. 
\eproof

\subsection{The case of $\kappa(K_{X}+B)=2$}
In this subsection, we focus on the case of $\kappa(K_{X}+B)=2$, which is the most difficult case.

Let $(X,B)$ be a projective lc threefold pair over an algebraically closed field  $k$ of characteristic $>3$ such that $K_{X}+B$ is nef and $\kappa(K_{X}+B)=2$. We pass to an uncountable base field.
After replacing $(X,B)$, we can assume that $(X,B)$ is $\Q$-factorial and dlt by Theorem \ref{dltmodification}.
Then one of the following cases holds:\\

\noindent Case \uppercase\expandafter{\romannumeral1}: $K_{X}+B-\varepsilon\lfloor B\rfloor$ is not pseudo-effective for any rational $\varepsilon>0$,\\

\noindent Case \uppercase\expandafter{\romannumeral2}: $K_{X}+B-\varepsilon\lfloor B\rfloor$ is pseudo-effective for any sufficiently small rational $\varepsilon>0$.\\

Note that by Lemma \ref{diagram}, $K_{X}+B$ is endowed with a map $h:X\to Z$ to a normal proper algebraic space $Z$ of dimension $2$. We will run several MMP which are $(K_{X}+B)$-trivial . It is clear that every step of such construction is still over $Z$.

\subsubsection{Proof of Case \uppercase\expandafter{\romannumeral1}}

In this part, we prove Case \uppercase\expandafter{\romannumeral1} (see Proposition \ref{dominatedcase}).
More precisely, we first prove that $\lfloor B\rfloor$ must dominate $Z$ in this case.
Then we deduce the semi-ampleness of $K_{X}+B$ by adjunction.

\blemma\label{findCartierdiv} Let $\phi:Z^{\prime}\to Z$ be a birational morphism from a $\Q$-factorial projective normal surface to a normal proper algebraic space of dimension $2$. Assume that $S$ is an effective Weil divisor on $Z^{\prime}$. Then we can take a $\Q$-Cartier $\Q$-divisor $A$ such that $A\geq S$ and $A\cdot E=0$ for any curve $E$ which is $\phi$-exceptional.
\elemma
\bproof We will write $A=S+H+\sum_{\alpha}a_{\alpha}C_{\alpha}$, where $H$ is a sufficiently ample effective divisor such that $S+H$ is ample, $C_{\alpha},\alpha\in I=\{ 1, 2, \ldots , r \}$ are all $\phi$-exceptional curves and $a_{\alpha}$ are some non-negative rational numbers. It is clear that $A\geq S$. We only need to choose appropriate $a_{\alpha}\geq 0$ such that $A\cdot E=0$ for any curve $E$ which is $\phi$-exceptional. 

Note that
\beq\begin{split}
  A\cdot C_{\beta}&=0, \beta \in I\\
\Longleftrightarrow\ \ \ \ \ \ \ \ \  (\sum_{\alpha}a_{\alpha}C_{\alpha})\cdot C_{\beta} &=-(S+H)\cdot C_{\beta}, \beta\in I\\
\Longleftrightarrow  [C_{\beta}\cdot C_{\alpha}]_{\alpha,\beta\in I}[a_{\alpha}]_{\alpha\in I}&=[-(S+H)\cdot C_{\beta}]_{\beta\in I},\notag
\end{split}\eeq
where $[C_{\beta}\cdot C_{\alpha}]_{\alpha,\beta\in I}$ is a matrix with element $C_{\beta}\cdot C_{\alpha}$ at row $\beta$ and column $\alpha$, and $[a_{\alpha}]_{\alpha\in I}, [-(S+H)\cdot C_{\beta}]_{\beta\in I}$ are column vectors with elements $a_{\alpha}, -(S+H)\cdot C_{\beta}$ at  rows $\alpha, \beta$, respectively.
Since $-(S+H)\cdot C_{\beta}<0$ for $\beta\in I$, to get a solution of $[a_{\alpha}]_{\alpha\in I}$ with $a_{\alpha}>0$ we only need to prove that the symmetric matrix $[C_{\beta}\cdot C_{\alpha}]_{\alpha,\beta\in I}$ is negative definite.

Consider a resolution of singularities $\phi^{\prime}:Z^{\prime\prime}\to Z^{\prime}$.
We first prove that the proposition holds for the morphism $\phi\circ\phi^{\prime}:Z^{\prime\prime}\to Z$. Let $C^{\prime}_{\alpha},\alpha\in J$ be all $\phi\circ\phi^{\prime}$-exceptional curves. Since $\phi\circ\phi^{\prime}$ is a contraction,
for any closed point $x\in Z$, $(\phi\circ\phi^{\prime})^{-1}(x)$ is connected. Hence different connected components of $\bigcup_{\alpha\in J}C^{\prime}_{\alpha}$ maps to different closed points. 
We apply \cite[Theorem 4.5]{artin1971algebraic} to the morphism $\phi\circ\phi^{\prime}$, then we know that the intersection matrix of any connected component of  $\bigcup_{\alpha\in J}C^{\prime}_{\alpha}$ is negative definite.
Note that the intersection matrix of $\bigcup_{\alpha\in J}C^{\prime}_{\alpha}$ is the direct sum of intersection matrices of all connected components of $\bigcup_{\alpha\in J}C^{\prime}_{\alpha}$. Hence the intersection matrix of $\bigcup_{\alpha\in J}C^{\prime}_{\alpha}$ is negative definite.

To prove that $[C_{\beta}\cdot C_{\alpha}]_{\alpha,\beta\in I}$ is negative definite, we only need to check $\phi^{\prime\ast}C_{\alpha},\alpha\in I$ are linearly independent. This is clear since we have $\phi^{\prime\ast}C_{\alpha}=\tilde{ C_{\alpha}}+E_{\alpha}$, where $\tilde{C_{\alpha}}$ are birational transforms of $C_{\alpha}$ and  $E_{\alpha}$ are $\phi^{\prime}$-exceptional $\Q$-divisors.
\eproof

\bproposition\label{terminationwhendegenerated} Let $(X,B)$ be a $\Q$-factorial projective dlt threefold pair over an algebraically closed field  $k$ of characteristic $>3$ with $\kappa(K_{X}+B)=2$. Assume $K_{X}+B$ is nef, and it is endowed with a map $h:X\to Z$. 
If $K_{X}+B-\varepsilon \lfloor B \rfloor$ is not pseudo-effective for any rational $\varepsilon>0$, then $\lfloor B\rfloor$ dominates $Z$.
\eproposition
\bproof 
We first prove the case when $X$ is terminal.
Since $K_{X}+B-\varepsilon \lfloor B \rfloor$ is not pseudo-effective for any rational $\varepsilon >0$, $K_{X}+(1-\varepsilon)B$ is not pseudo-effective for any rational $\varepsilon >0$. Then by Definition \ref{MMPKtrivial} we can run a $K_{X}$-MMP which is $(K_{X}+B)$-trivial. By Lemma \ref{terminaltermination} it terminates with a pair $(X^{\prime},B^{\prime})$ since $X$ is terminal. 
Moreover, since $\kappa(K_{X}+(1-\varepsilon)B)=\kappa(K_{X^{\prime}}+(1-\varepsilon)B^{\prime})$ and $K_{X}+(1-\varepsilon)B$ is not pseudo-effective for any small rational $\varepsilon >0$, $K_{X^{\prime}}+(1-\varepsilon)B^{\prime}$ is not nef for any small rational $\varepsilon>0$ by Theorem \ref{nonvan}.
Hence this MMP terminates 
with a Mori fibre space  
\beq\begin{tikzcd}
X\arrow[r,dashrightarrow,"f"]& Y\arrow[d,"h^{\prime}"]\\
&Z^{\prime}\notag
\end{tikzcd}\eeq 
Denote the birational transform of $B$ on $Y$ by $B_{Y}$.
Note that $K_{Y}+B_{Y}$ is endowed with a map $h_{Y}:Y\to Z$ and $h_{Y}$ factors through $h^{\prime}$ since $h^{\prime}$ only contracts curves which are $(K_{Y}+B_{Y})$-trivial. In other words, we have a commutative diagram 
\beq\begin{tikzcd}
&Y\arrow[dr,"h_{Y}"]\arrow[dl,swap,"h^{\prime}"]& \\
Z^{\prime}\arrow[rr,"\phi "]& & Z.\notag
\end{tikzcd}\eeq
Note that $h^{\prime}$ is equidimensional, $Z^{\prime}$ is $\Q$-factorial and $\phi$ is   a birational map. 
Applying Lemma \ref{findCartierdiv} to $\phi$, we get a $\Q$-Cartier $\Q$-divisor $A$ on $Z^{\prime}$ such that $A\geq h^{\prime}(\lfloor B_{Y}\rfloor)$ and $h^{\prime}(F)\cdot A=0$ for any $h_{Y}$-exceptional divisor $F$.
Note that
$$\kappa(K_{Y}+B_{Y}-\lfloor B_{Y}\rfloor+a h^{\prime\ast}A)\geq \kappa(K_{Y}+B_{Y})= \kappa(K_{X}+B)=2$$
for some integer $a>0$. Hence there exists an effective $\Q$-divisor 
$$M\sim_{\Q}K_{Y}+B_{Y}-\lfloor B_{Y}\rfloor+a h^{\prime\ast}A$$
 such that $M\cdot C=(K_{Y}+B_{Y}-\lfloor B_{Y}\rfloor)\cdot C$ for any curve $C$ in the fibre of $h_{Y}$.
In other words, flips of a $(K_{Y}+B_{Y}-\lfloor B_{Y}\rfloor)$-MMP which is $(K_{Y}+B_{Y})$-trivial are all $M$-filps.
Therefore, by Theorem \ref{Waldrontermination} a $(K_{Y}+B_{Y}-\lfloor B_{Y}\rfloor)$-MMP which is $(K_{Y}+B_{Y})$-trivial terminates with a Mori fibre space
\beq\begin{tikzcd}
Y\arrow[r,dashrightarrow,"f^{\prime}"]& Y^{\prime}\arrow[d,"h^{\prime\prime}"]\\
&Z^{\prime\prime}\notag
\end{tikzcd}\eeq
such that $f^{\prime}_{\ast}(\lfloor B_{Y}\rfloor)$ dominates $Z^{\prime\prime}$
by Lemma \ref{Output}.
Note that $K_{Y^{\prime}}+B_{Y^{\prime}}$ is endowed with a map $h_{Y^{\prime}}:Y^{\prime}\to Z$ and $h_{Y^{\prime}}$ factors through $h^{\prime\prime}$ since $h^{\prime\prime}$ only contracts curves which are $(K_{Y^{\prime}}+B_{Y^{\prime}})$-trivial.
Therefore $f^{\prime}_{\ast}(\lfloor B_{Y}\rfloor)$ , and hence $\lfloor B \rfloor$ dominate $Z$.

Now we turn to the general case. By Theorem \ref{dltmodification} we can take a  
dlt modification $\mu:(X^{\prime\prime},B^{\prime\prime})\to (X,B)$ such that $(X^{\prime\prime},B^{\prime\prime})$ is $\Q$-factorial and dlt, and $X^{\prime\prime}$ is terminal. 
If $K_{X}+B-\varepsilon \lfloor B \rfloor$ is not pseudo-effective for any rational $\varepsilon >0$, then $K_{X^{\prime\prime}}+B^{\prime\prime}-\varepsilon\lfloor B^{\prime\prime}\rfloor$ is not pseudo-effective for any rational $\varepsilon >0$, since 
$$\mu_{\ast}(K_{X^{\prime\prime}}+B^{\prime\prime}-\varepsilon\lfloor B^{\prime\prime}\rfloor)=K_{X}+B-\varepsilon \lfloor B \rfloor.$$
By the last paragraph, $\lfloor B^{\prime\prime}\rfloor$ dominates $Z$.
Note that $\lfloor B^{\prime\prime}\rfloor$ dominates $Z$ if and only if $\lfloor B\rfloor$ dominates $Z$ since $Z$ is of dimension $2$ and $\mu$ is an isomorphism over a big open subset of $X$.
Hence we have $\lfloor B\rfloor$ dominates $Z$.
\eproof

Now we can prove Case \uppercase\expandafter{\romannumeral1}.

\bproposition\label{dominatedcase} Let $(X,B)$ be a $\Q$-factorial projective dlt threefold pair over an algebraically closed field  $k$ of characteristic $>3$ such that $K_{X}+B$ is nef and $\kappa(K_{X}+B)=2$.  
If $K_{X}+B-\varepsilon \lfloor B \rfloor$ is not pseudo-effective for any rational $\varepsilon>0$, then $K_{X}+B$ is semi-ample. 
\eproposition
\bproof We pass to an uncountable base field. By Lemma \ref{diagram}, $K_{X}+B$ is endowed with a map $h:X\to Z$ to a normal proper algebraic space $Z$ of dimension $2$.
Now by Proposition \ref{terminationwhendegenerated}, $\lfloor B\rfloor$ dominates $Z$.

Since $(K_{X}+B)|_{G}\equiv 0$, where $G$ is the generic fibre of $h$ and $G$ is of dimension $1$, we have $(K_{X}+B)|_{G}\sim_{\Q} 0$ by the abundance for curves. 
Then we can apply Lemma \ref{diagram} to get a commutative diagram
\beq\begin{tikzcd}
X_{1}\arrow[r,"\phi "]\arrow[d,swap,"h_{1}"]&X\arrow[d,"h"] \\
Z_{1}\arrow[r,"\psi"] &Z\notag
\end{tikzcd}
\eeq
where $Z_{1}$ is a smooth projective surface, $X_{1}$ is a normal projective threefold, $\phi,\psi$ are birational morphisms and $h_{1}: X_{1}\to Z_{1}$ is a fibration. Moreover, there exists a nef and big $\Q$-divisor $D_{1}$ on $Z_{1}$ such that $\phi^{\ast}(K_{X}+B)\sim_{\Q} h^{\ast}_{1}D_{1}$. To show $K_{X}+B$ is semi-ample, it suffices to show $D_{1}$ is semi-ample. 

Let $B_{1}$ be the birational transform of $B$ on $X_{1}$. Since $\lfloor B\rfloor$ dominates $Z$, we have $\lfloor B_{1}\rfloor$ dominates $Z_{1}$. 
Moreover we have $\phi^{\ast}(K_{X}+B)|_{\lfloor B_{1}\rfloor}$ is semi-ample since $(K_{X}+B)|_{\lfloor B\rfloor}$ is semi-ample by Theorem \ref{adjunction}.
Hence by Lemma \ref{pullbacksemi-ample}, $D_{1}$, and hence $K_{X}+B$ are semi-ample.
\eproof

\subsubsection{Proof of Case \uppercase\expandafter{\romannumeral2}}

In this part, we prove Case \uppercase\expandafter{\romannumeral2} (see Proposition \ref{Case2}).
First, we prove this case when $K_{X}+B$ is endowed with an equidimensional map $h:X\to Z$. For the general case, we modify the pair $(X,B)$ by running several MMP which are $(K_{X}+B)$-trivial so that all $h$-exceptional prime divisors are connected components of $\lfloor B\rfloor$.
Then after further modification we can construct an equidimensional fibration
$h_{\varepsilon}:X\to Z_{\varepsilon}$ to a normal projective surface.
Finally, we descend $K_{X}+B$ to $Z_{\varepsilon}$ and prove its semi-ampleness.

\bproposition\label{mainthmwhenequidim} Let $D$ be a nef $\Q$-divisor on $X$ with $\kappa(X,D)=2$, where $X$ is a $\Q$-factorial normal projective threefold over an uncountable algebraically closed field $k$ of characteristic $>0$. Suppose that $D$ is endowed with an equidimensional map $h:X\to Z$ such that $D|_{G}\sim_{\Q}0$, where $G$ is the generic fibre of $h$.
Then $Z$ is a projective variety and $D$ is semi-ample.
\eproposition
\bproof By Lemma \ref{diagram},
there is a commutative diagram as following
\beq\begin{tikzcd}
X_{1}\arrow[r,"\phi "]\arrow[d,swap,"h_{1}"]&X\arrow[d,"h"] \\
Z_{1}\arrow[r,"\psi"] &Z\notag
\end{tikzcd}
\eeq
where $Z_{1}$ is a smooth projective surface, $X_{1}$ is a normal projective threefold, $\phi,\psi$ are birational morphisms and $h_{1}: X_{1}\to Z_{1}$ is an equidimensional fibration. Moreover, there exists a nef and big $\Q$-divisor $D_{1}$ on $Z_{1}$ such that $\phi^{\ast}D\sim_{\Q} h^{\ast}_{1}D_{1}$.

Since $Z$ is a normal proper algebraic space of dimension $2$, there exists an open set $U\subseteq Z$ such that $U$ is a smooth quasi-projective variety and $T:=Z\backslash U$ consists of finitely many closed points on $Z$. By Lemma \ref{descend} we have $D|_{h^{-1}(U)}$ is $\Q$-linearly trivial over $U$ since $h$ is equidimensional and $D|_{G}\sim_{\Q}0$.
Now we take a very ample divisor $S$ on $X$, which does not contain any component of $h^{-1}(T)$. Then we have the following commutative diagram 
\beq\begin{tikzcd}
S_{1}^{\nu}\arrow[rr,"\mathrm{normalization} "]\arrow[d,"\phi_{S^{\nu}}"]& &S_{1}=\phi^{-1}S\arrow[r]\arrow[d,"\phi_{S}"]&X_{1}\arrow[r,"h_{1}"]\arrow[d,"\phi "]&Z_{1} \arrow[d,"\psi"]\\
S^{\nu}\arrow[rr,"\mathrm{normalization}"]& &S\arrow[r]&X\arrow[r,"h"] &Z.              \notag
\end{tikzcd}
\eeq
 The $\Q$-divisor $D|_{S^{\nu}}$ is nef and big. Consider the exceptional locus $\E(D|_{S^{\nu}})$. It is, the union of finitely many $D$-numerically trivial curves on $S^{\nu}$. Note that $S\cap h^{-1}(T)$ contains no curve by our construction. Hence the image of $\E(D|_{S^{\nu}})$, via the natural map $S^{\nu} \to X$, is contained in finitely many fibers of $h$ over some closed points in $U$. Therefore $(D|_{S^{\nu}})|_{\E(D|_{S^{\nu}})}$ is semi-ample, and by Theorem \ref{Keelsemi-ample} $D|_{S^{\nu}}$ is semi-ample. 

Denote the natural map $S_{1}^{\nu}\to Z_{1}$ by $\sigma$.
Since $D|_{S^{\nu}}$ is semi-ample, we know that
$$\phi_{S^{\nu}}^{\ast} D|_{S^{\nu}}\sim_{\Q}\sigma^{\ast}D_{1}$$
is semi-ample. Then by Lemma \ref{pullbacksemi-ample} we have $D_{1}$ is semi-ample.
Hence $\phi^{\ast}D\sim_{\Q} h^{\ast}_{1}D_{1}$ is semi-ample. Again by Lemma \ref{pullbacksemi-ample} it follows that $D$ is semi-ample.
Moreover, $D$ induces the morphism $h:X\to Z$. Hence $Z$ is projective.
\eproof

This proposition proves Case \uppercase\expandafter{\romannumeral2} when $K_{X}+B$ is endowed with an equidimensional map $h:X\to Z$ by letting $D=K_{X}+B$. In general, this equidimensionality condition may fail. We need to modify the pair $(X,B)$. To do this, we need the following lemmas.

\blemma\label{Fnotnef} Let $D$ be a nef $\Q$-divisor on $X$ with $\kappa(X,D)=2$, where $X$ is a $\Q$-factorial normal projective threefold over an uncountable algebraically closed field $k$ of characteristic $>0$. Suppose that $D$ is endowed with a map $h:X\to Z$ such that $D|_{G}\sim_{\Q}0$, where $G$ is the generic fibre of $h$.
Then any $h$-exceptional prime divisor $F$ is not nef.
\elemma
\bproof By Lemma \ref{diagram}, we have the following commutative diagram
\beq\begin{tikzcd}
X_{1}\arrow[r,"\phi "]\arrow[d,swap,"h_{1}"]&X\arrow[d,"h"] \\
Z_{1}\arrow[r,"\psi"] &Z \notag
\end{tikzcd}
\eeq
where $Z_{1}$ is a smooth projective surface, $X_{1}$ is a normal projective threefold, $\phi,\psi$ are birational morphisms and $h_{1}: X_{1}\to Z_{1}$ is an equidimensional fibration such that, there exists a nef and big $\Q$-divisor $D_{1}$ on $Z_{1}$ such that $\phi^{\ast}D\sim_{\Q} h^{\ast}_{1}D_{1}$.

First by the definition of EWM we have $D$ is numerically trivial on $F$. 
Let $F_{1}$ be the birational transform of $F$ on $X_{1}$.
Since $D_{1}$ is a nef and big $\Q$-divisor on $Z_{1}$, we can write $D_{1}\sim_{\Q}A+E_{1}$ such that $A$ is an ample effective $\Q$-divisor, and $E_{1}$ is an effective $\Q$-divisor. Moreover, we can choose $A$ such that $\mathrm{Supp}(h_{1}^{\ast}A)$ doesn't contain any component of $\mathrm{Supp}(\phi^{\ast}F)\cup \mathrm{Exc}(\phi)$  since $A$ is ample.
We take a $\Q$-effective divisor $\Delta$ such that $D\sim_{\Q}\Delta$ and $\phi^{\ast}\Delta=h_{1}^{\ast}(A+E_{1})$. 
 
Now we take a very ample divisor $H_{1}$ on $X_{1}$.
Since $h_{1}^{\ast}A\cdot F_{1}\cdot H_{1}>0$, we have $\mathrm{Supp}(h_{1}^{\ast}A)\cap F_{1}\neq \emptyset$.
Let $A_{X}$ be the birational transform of $\mathrm{Supp}(h_{1}^{\ast}A)$ on $X$. Then its intersection with $F$ is of dimension one by our choice of $A$. If we take a very ample divisor $H$ on $X$, it is clear that $A_{X}\cdot F\cdot H>0$. 
Note that $\Delta\cdot F\cdot H=0$ and $A_{X}\subseteq \mathrm{Supp}\ \Delta$.
It implies that $F\subseteq\mathrm{Supp}\ \Delta$ and $F\cdot F \cdot H<0$.
\eproof

\blemma\label{perturb} Let $(X,B)$ be a $\Q$-factorial projective lc threefold pair over an algebraically closed field  $k$ of characteristic $>3$, and $D$ be an effective $\Q$-divisor such that $\mathrm{Supp}\ D\subseteq \mathrm{Supp}\ B$.
 Assume that $K_{X}+B$ is nef and $K_{X}+B-\varepsilon D$ is pseudo-effective for any sufficiently small rational $\varepsilon>0$. 
 Then we have
 
\noindent (1) $\kappa(K_{X}+B-\varepsilon D)=\kappa(K_{X}+B)$ for any sufficiently small rational $\varepsilon>0$,

\noindent (2) if $D\subseteq \lfloor B\rfloor$ is a reduced divisor, then any $(K_{X}+B- D)$-MMP which is $(K_{X}+B)$-trivial terminates with a pair $(X^{\prime},B^{\prime})$ such that $K_{X^{\prime}}+B^{\prime}-\varepsilon D^{\prime}$ is nef for any sufficiently small rational $\varepsilon>0$, where $D^{\prime}$ is the birational transform of $D$ on $X^{\prime}$,

\noindent (3) if $D\subseteq \lfloor B\rfloor$ is a prime divisor, then $D$ is not contracted by any $(K_{X}+B- D)$-MMP which is $(K_{X}+B)$-trivial.
\elemma
\bproof 
Since $K_{X}+B-\varepsilon D$ is pseudo-effective for any sufficiently small rational $\varepsilon>0$, by Theorem \ref{nonvan} we have $K_{X}+B-\varepsilon D$ is effective for any sufficiently small rational $\varepsilon>0$. 
Hence there exists an effective $\Q$-divisor 
$\Delta_{\varepsilon}\sim_{\Q}K_{X}+B-2\varepsilon D$ for a sufficiently small rational $\varepsilon>0$.
Then we have
$$K_{X}+B\sim_{\Q}\Delta_{\varepsilon}+2\varepsilon D, K_{X}+B-\varepsilon D \sim_{\Q}\Delta_{\varepsilon}+\varepsilon D.$$ 
This proves (1) since effective divisors with the same support have the same Kodaira dimension. 

Assume that $D\subseteq \lfloor B\rfloor$ is a reduced divisor. Note that for any sufficiently small rational $\varepsilon>0$, $K_{X}+B-\varepsilon D$ is pseudo-effective and every step of a $(K_{X}+B-\varepsilon D)$-MMP which is $(K_{X}+B)$-trivial
is a step of a $(K_{X}+B-\varepsilon D)$-MMP. We choose a sufficiently small rational $\varepsilon_{0}>0$.
By Theorem \ref{termination}, we have 
a $(K_{X}+B-\varepsilon_{0} D)$-MMP which is $(K_{X}+B)$-trivial terminates with 
a pair $(X^{\prime},B^{\prime})$ such that $K_{X^{\prime}}+B^{\prime}-\varepsilon D^{\prime}$ is nef for any sufficiently small rational $\varepsilon>0$, where $D^{\prime}$ is the birational transform of $D$ on $X^{\prime}$.
Since any $(K_{X}+B- D)$-MMP which is $(K_{X}+B)$-trivial is a $(K_{X}+B-\varepsilon_{0} D)$-MMP which is $(K_{X}+B)$-trivial, we have (2) holds.

Assume moreover that $D$ is a prime divisor. By (2), a $(K_{X}+B- D)$-MMP which is $(K_{X}+B)$-trivial terminates with a pair $(X^{\prime},B^{\prime})$ such that $K_{X^{\prime}}+B^{\prime}-\varepsilon D^{\prime}$ is nef for any sufficiently small rational $\varepsilon>0$, where $D^{\prime}$ is the birational transform of $D$ on $X^{\prime}$.
We take a common resolution of $X$ and $X^{\prime}$ 
\beq\begin{tikzcd}
&W\arrow[rd,"\phi_{2} "]\arrow[ld,swap,"\phi_{1}"]& \\
X\arrow[rr,dashrightarrow,"f"]& &X^{\prime} \notag
\end{tikzcd}
\eeq
Note that since every step of a $(K_{X}+B- D)$-MMP which is $(K_{X}+B)$-trivial is a step  of a $(K_{X}+B- D)$-MMP, we have
$$\phi_{1}^{\ast}(K_{X}+B-D)\sim_{\Q}\phi^{\ast}_{2}(K_{X^{\prime}}+B^{\prime}-D^{\prime})+E,$$
where $E$ is an effective $\phi_{2}$-exceptional $\Q$-divisor. 
It implies that
$$-\phi^{\ast}_{1}D-E\sim_{\Q}\phi^{\ast}_{2}(K_{X^{\prime}}+B^{\prime}-D^{\prime})-\phi^{\ast}_{1}(K_{X}+B).$$
Applying the negativity lemma (\cite[Lemma 3.39]{kollar1998birational}) to $\phi_{2}$, we know that
$$-\phi_{2\ast}\phi^{\ast}_{1}D \neq 0.$$
Hence $D$ is not contracted by $f$, i.e. (3) holds.
\eproof

Now we can prove Case \uppercase\expandafter{\romannumeral2}.

\bproposition\label{Case2} Let $(X,B)$ be a $\Q$-factorial projective dlt threefold pair over an algebraically closed field  $k$ of characteristic $>3$ such that $K_{X}+B$ is nef and $\kappa(K_{X}+B)=2$.  
If $K_{X}+B-\varepsilon \lfloor B \rfloor$ is pseudo-effective for any sufficiently small rational $\varepsilon>0$, then $K_{X}+B$ is semi-ample.
\eproposition
\bproof We pass to an uncountable base field. By Proposition \ref{dominatedcase},
$K_{X}+B$ is endowed with a map $h:X\to Z$ to an algebraic space $Z$ of dimension $2$.

\noindent Step 1. We contract all $h$-exceptional prime divisors which have empty intersection with $\lfloor B\rfloor$.

Let $F$ be a $h$-exceptional prime divisor such that $F\cap \lfloor B\rfloor=\emptyset$, then we can choose a sufficiently small rational $\varepsilon$ such that $(X,B+\varepsilon F)$ is still dlt. Note that by Lemma \ref{Fnotnef} we have $K_{X}+B+\varepsilon F$ is not nef since $K_{X}+B$ is numerically trivial on $F$. We run a $(K_{X}+B+\varepsilon F)$-MMP as follows.

For the first step, the extremal ray is $(K_{X}+B)$-numerically trivial since any curve which is $(K_{X}+B+\varepsilon F)$-negative must be contained in $F$. If it is a divisorial contraction, then $F$ is contracted and the process terminates. Otherwise, we get a flip
$$\mu:(X,B+\varepsilon F) \ratmap (X^{+},B^{+}+\varepsilon F^{+})$$
 such that $F^{+}\neq 0$. Note that $K_{X^{+}}+B^{+}+ F^{+}$ is still not nef.
By Theorem \ref{termination} the process must terminate, hence $F$ is contracted after finitely many steps. Since at every step we only contract $(K_{X}+B)$-trivial curves, we can replace $(X,B)$ by the output of this process. Moreover, since the number of $h$-exceptional prime divisors is finite, we can repeat this process until every $h$-exceptional divisor intersects $\lfloor B \rfloor$.

From now on, we can assume that every $h$-exceptional divisor intersects $\lfloor B\rfloor$. \\

\noindent Step 2. We reduce the proposition to the case when all $h$-exceptional prime divisors are connected components of $\lfloor B\rfloor$.

To this end, let $S\subseteq \lfloor B\rfloor$ be a prime divisor such that there exists a $h$-exceptional divisor $F$ whose intersection with $S$ is of dimension one. By Definition \ref{MMPKtrivial} we run a $(K_{X}+B-S)$-MMP which is $(K_{X}+B)$-trivial.
By Lemma \ref{perturb}, it terminates with a pair $(X_{1},B_{1})$ such that $K_{X_{1}}+B_{1}-\varepsilon S_{1}$ is nef for any sufficiently small rational $\varepsilon>0$, where $S_{1}$ is the birational transform of $S$ on $X_{1}$.
Moreover, $S_{1}\neq 0$.

After replacing $(X,B), S$ by $(X_{1},B_{1}), S_{1}$ ($(X,B)$ may no longer be dlt) , we can assume that $K_{X}+B-\varepsilon S$ is nef for any sufficiently small rational $\varepsilon>0$. Let $F$ be a $h$-exceptional prime divisor such that it has non-empty intersection with $S$. Since $K_{X}+B$ is numerically-trivial on $F$, we have $-S$ is nef on $F$, which implies that $S=F$. It is to say that after this process, there is no $h$-exceptional divisor $F$ whose intersection with $S$ is of dimension one.

Since the number of $h$-exceptional prime divisors is finite and it decreases strictly under the above process, we can repeat this process until there is no prime divisor $S\subseteq \lfloor B\rfloor$ such that there exists a $h$-exceptional divisor $F$ whose intersection with $S$ is of dimension one.

From now on, we can assume that all $h$-exceptional prime divisors are connected components of $\lfloor B\rfloor $.\\

\noindent Step 3. We further modify $(X,B)$ and construct an equidimensional fibration $h_{\varepsilon}:X\to Z_{\varepsilon}$.

First, let $F_{h}$ be the reduced $h$-exceptional divisor and run a $(K_{X}+B-F_{h})$-MMP which is $(K_{X}+B)$-trivial by Definition \ref{MMPKtrivial}. After replacing $(X,B)$ by the output of this process, we can assume that $K_{X}+B-\varepsilon F_{h}$ is nef for any sufficiently small rational $\varepsilon>0$ as at Step 2.

We choose a sufficiently small rational $\varepsilon>0$.
Note that by Lemma \ref{perturb} we have $\kappa(K_{X}+B-\varepsilon F_{h})=\kappa(K_{X}+B)=2$.
Hence by Lemma \ref{diagram}, $K_{X}+B-\varepsilon F_{h}$ is endowed with a map $h_{\varepsilon}:X\to Z_{\varepsilon}$.
We claim that there exists a commutative diagram
\beq\begin{tikzcd}
X\arrow[dr,"h_{\varepsilon} "]\arrow[rr,"h"]& &Z \\
 &Z_{\varepsilon}\arrow[ur,"\psi_{\varepsilon} "] & \notag
\end{tikzcd}
\eeq
We only need to prove that any curve contracted by $h_{\varepsilon}$ is contracted by $h$. Let $C_{1}$ be a curve contracted by $h_{\varepsilon}$, i.e. $(K_{X}+B-\varepsilon F_{h})\cdot C_{1}=0$. 

If $C_{1}\cap F_{h}=\emptyset$, then $(K_{X}+B-\varepsilon F_{h})\cdot C_{1}=0$ implies $(K_{X}+B)\cdot C_{1}=0$. Hence $C_{1}$ is contracted by $h$.
If $C_{1}\cap F_{h}\neq\emptyset$ and $C_{1}\not\subseteq F_{h}$, then we have $C_{1}\cdot F_{h}>0$ . But $K_{X}+B-2\varepsilon F_{h}$ is nef as well, i.e. 
$$(K_{X}+B-2\varepsilon F_{h})\cdot C_{1}=-\varepsilon F_{h}\cdot C_{1}\geq 0.$$ We obtain a contradiction. Finally, if $C_{1}\subseteq F_{h}$, then  $C_{1}$ is always contracted by $h$. 

We prove that $h_{\varepsilon}$ is actually equidimensional. By the above diagram we know that exceptional divisors of $h_{\varepsilon}$ have to be exceptional divisors of $h$. Hence all $h_{\varepsilon}$-exceptional divisors are supported in $F_{h}$. If $F$ is a prime $h_{\varepsilon}$-exceptional divisor , we have both $K_{X}+B$ and $K_{X}+B-\varepsilon F_{h}$ are numerically trivial on $F$ , and hence $F_{h}$ is numerically trivial on $F$, which is impossible since $F$ is not nef by Lemma \ref{Fnotnef} and $F$ is a connected component of $F_{h}$.\\

\noindent Step 4. Descend $K_{X}+B$ to $Z_{\varepsilon}$ and prove its semi-ampleness.

By Proposition \ref{mainthmwhenequidim}, we have $K_{X}+B-\varepsilon F_{h}$ is semi-ample and $Z_{\varepsilon}$ is a projective variety. Moreover by Lemma \ref{descend} $K_{X}+B$ descends to a nef and big divisor $D_{\varepsilon}$ on $Z_{\varepsilon}$ since $h_{\varepsilon}$ is equidimensional and $Z_{\varepsilon}$ is $\Q$-factorial by \cite[Proposition 3.3]{waldron2017finite}.

By the projection formula
for any curve $\Gamma\subseteq \E(D_{\varepsilon})$ we have $K_{X}+B$ is numerically trivial on $h_{\varepsilon}^{-1}(\Gamma)$. However by our assumption, $h_{\varepsilon}^{-1}(\Gamma)$ has to be contained in $F_{h}$.
Hence it is clear that
$$\E(D_{\varepsilon})\subseteq h_{\varepsilon}(F_{h}).$$
Since $h_{\varepsilon}$ is equidimensional, we have $h_{\varepsilon}^{-1}(h_{\varepsilon}(F_{h}))$ is the union of finitely many prime divisors. All these prime divisors are exceptional divisors of $h$ since $\psi_{\varepsilon}\circ h_{\varepsilon}(F_{h})$ is of dimension $0$.  Hence, we have 
$$h_{\varepsilon}^{-1}(\E(D_{\varepsilon}))\subseteq h_{\varepsilon}^{-1}(h_{\varepsilon}(F_{h}))= F_{h}.$$

We take a dlt modification $g: (X^{\prime},B^{\prime})\to (X,B)$ such that $g$ only extracts  prime divisors $E$ with discrepancies $a(E,X,B)=-1$ by \cite[Lemma 7.7]{birkar2016existence} and \cite{hacon2022minimal}.  
Then we have
$$(h_{\varepsilon}\circ g)^{-1}(\E(D_{\varepsilon}))\subseteq g^{-1}(F_{h})\subseteq g^{-1}(\lfloor B\rfloor)\subseteq \lfloor B^{\prime}\rfloor.$$
Since $(K_{X^{\prime}}+B^{\prime})|_{\lfloor B^{\prime}\rfloor}$ is semi-ample by Theorem \ref{adjunction}, we have 
$$(K_{X^{\prime}}+B^{\prime})|_{(h_{\varepsilon}\circ g)^{-1}(\E(D_{\varepsilon}))}$$
is semi-ample. Then by \cite[Lemma 7.1]{birkar2017existence} we have $D_{\varepsilon}$, and hence $K_{X}+B$ are semi-ample.
\eproof

\subsection{Proof of Theorem \ref{mainthm}}

\bproof Case of $\kappa(X,K_{X}+B)=3$: In this case $K_{X}+B$ is nef and big, hence the proposition holds by \cite[Theorem 1.1]{waldron2018lmmp} and \cite{hacon2022minimal}.\\

\noindent Case of $\kappa(X,K_{X}+B)=2$:
After replacing $(X,B)$ by its dlt modification, we can assume that $(X,B)$ is a $\Q$-factorial dlt pair by Theorem \ref{dltmodification}.
Then the proposition follows from Proposition \ref{dominatedcase} and Proposition \ref{Case2}.\\

\noindent Case of $\kappa(K_{X}+B)=1$: The proof is similar to the case of $\kappa(K_{X}+B)=2$ but easier.
 
After replacing $(X,B)$ by its dlt modification, we can assume that $(X,B)$ is a $\Q$-factorial dlt pair and $X$ is terminal by Theorem \ref{dltmodification}. 
Then we have either

\noindent (1): $K_{X}+B-\varepsilon\lfloor B\rfloor$ is not pseudo-effective for any rational $\varepsilon>0$, or

\noindent (2): $K_{X}+B-\varepsilon\lfloor B\rfloor$ is pseudo-effective for any sufficiently small rational $\varepsilon>0$.

In the case of (1), since $K_{X}+B-\varepsilon \lfloor B \rfloor$ is not pseudo-effective for any rational $\varepsilon >0$, $K_{X}+(1-\varepsilon)B$ is not pseudo-effective for any rational $\varepsilon >0$. Then by Definition \ref{MMPKtrivial} we can run a $K_{X}$-MMP which is $(K_{X}+B)$-trivial. By Lemma \ref{terminaltermination} it terminates with a pair $(X^{\prime},B^{\prime})$ since $X$ is terminal. 
Moreover, since 
$$\kappa(K_{X}+(1-\varepsilon)B)=\kappa(K_{X^{\prime}}+(1-\varepsilon)B^{\prime})$$
and $K_{X}+(1-\varepsilon)B$ is not pseudo-effective for any small rational $\varepsilon >0$, $K_{X^{\prime}}+(1-\varepsilon)B^{\prime}$ is not nef for any small rational $\varepsilon>0$ by Theorem \ref{nonvan}.
Hence this $K_{X}$-MMP which is $(K_{X}+B)$-trivial terminates with a Mori fibre space.
Then we have $n(K_{X}+B)\leq 2$. By Proposition \ref{weakcan}, 
 $n(K_{X}+B)=1$ since $\kappa(K_{X}+B)=1$. Then a nef reduction map of $K_{X}+B$, which exists by Theorem \ref{nefreductionmap},  is an equidimensional fibration to a normal curve.
Hence we can descend $K_{X}+B$ to an ample divisor on the curve by Lemma \ref{descend}.

In the case of (2), by Definition \ref{MMPKtrivial} we run a $(K_{X}+B-\lfloor B\rfloor)$-MMP which is $(K_{X}+B)$-trivial which terminates by Lemma \ref{perturb}, and replace $(X,B)$ by the output. $(X,B)$ may no longer be dlt and $X$ may no longer be terminal. However, we can assume that $K_{X}+B-\varepsilon\lfloor B\rfloor$ is nef and $(X,B-\varepsilon\lfloor B\rfloor)$ is klt for any sufficiently small rational $\varepsilon>0$.
By Lemma \ref{perturb} we have 
$$\kappa(K_{X}+B-\varepsilon\lfloor B\rfloor)=\kappa(K_{X}+B)=1$$ 
for any sufficiently small rational $\varepsilon> 0$. 
We choose a sufficiently small rational $\varepsilon>0$ such that $K_{X}+B-2\varepsilon\lfloor B\rfloor$ is nef and $\kappa(K_{X}+B-\varepsilon\lfloor B\rfloor)=1$.
Then by Theorem \ref{known}, $|m(K_{X}+B-\varepsilon\lfloor B\rfloor)|$ induces a fibration $h^{\prime}:X\to Z^{\prime}$ for a sufficiently divisible positive integer $m$ since $(X,B-\varepsilon\lfloor B\rfloor)$ is klt. 
Denote the generic fibre of $h^{\prime}$ by $G$. 
By Theorem \ref{surfaceabundance}, $(K_{X}+B-2\varepsilon\lfloor B\rfloor)|_{G}$ is semi-ample. 
Note that $(K_{X}+B-\varepsilon\lfloor B\rfloor)|_{G}\sim_{\Q}0$. It implies that $(K_{X}+B-2\varepsilon\lfloor B\rfloor)|_{G}\sim_{\Q}0$, and hence $(K_{X}+B)|_{G}\sim_{\Q}0$. 
Then by Lemma \ref{descend}, $K_{X}+B$ descends to an ample divisor on $Z^{\prime}$. Hence $K_{X}+B$ is semi-ample.
\eproof

\section{Applications}

In this section, we complete the proofs of the statements in the introduction.

\btheorem\label{finitegeneration} Let $(X, B)$ be a projective lc threefold pair over an algebraically closed field  $k$ of characteristic $>3$. Then the log canonical ring
$$R(K_{X}+B)=\oplus_{m=0}^{\infty} H^{0}(\lfloor m(K_{X}+B)\rfloor )$$
is finitely generated.
\etheorem
\bproof If $\kappa(K_{X}+B)=0$ or $-\infty$, the assertion is trivial. Otherwise, we have $\kappa(K_{X}+B)\geq 1$. After replacing $(X,B)$ by its log minimal model by Theorem \ref{lcmmp}, we can assume that $K_{X}+B$ is nef. Then the assertion follows from Theorem \ref{mainthm}.
\eproof

\btheorem\label{nefdim} Let $(X,B)$ be a projective lc threefold pair over an algebraically closed field  $k$ of characteristic $>3$. If $K_{X}+B$ is nef and $n(X,K_{X}+B)\leq 2$, then $K_{X}+B$ is semi-ample.
\etheorem
\bproof Case of $n(K_{X}+B)=0$: By Theorem \ref{nonvan} we have $\kappa(K_{X}+B)\geq 0$. Hence we have
$$\kappa(K_{X}+B)=n(K_{X}+B)=0.$$
Therefore $K_{X}+B\sim_{\Q}0$.\\

\noindent Case of $n(K_{X}+B)=1$: Let $\phi:X\dashrightarrow Z$ be a nef reduction map, which exists by Theorem \ref{nefreductionmap},. Since $Z$ is a normal curve, $X$ is normal and $\phi$ is proper over the generic point $\mu$ of $Z$, we have $\phi$ is indeed a morphism. Note that $(K_{X}+B)|_{G}\sim_{\Q}0$ by Theorem \ref{surfaceabundance}, where $G$ is the generic fibre of $\phi$. Since $\phi$ is equidimensional, we have $K_{X}+B\sim_{\Q}f^{\ast}A$ for an ample divisor on $Z$ by Lemma \ref{descend}. Hence $K_{X}+B$ is semi-ample.\\

\noindent Case of $n(K_{X}+B)=2$: By Proposition \ref{weakcan}, we have $\kappa(K_{X}+B)=2$. Then the proposition follows from Theorem \ref{mainthm}.
\eproof

\btheorem\label{Alb}  Let $(X,B)$ be a projective lc threefold pair over an algebraically closed field $k$ of characteristic $>3$.  If $K_{X}+B$ is nef and $\mathrm{dim}\ \mathrm{Alb}(X)\neq 0$, then $K_{X}+B$ is semi-ample.
\etheorem
\bproof After replacing $(X,B)$ by its dlt modification, we can assume that $(X,B)$ is a $\Q$-factorial dlt pair and $X$ is terminal by Theorem \ref{dltmodification}.
Moreover, by Theorem \ref{nonvan} and Theorem \ref{mainthm} we can assume that $\kappa(K_{X}+B)=0$.
By Definition \ref{MMPKtrivial} we run a $K_{X}$-MMP which is $(K_{X}+B)$-trivial, which terminates by Lemma \ref{terminaltermination} since $X$ is terminal. 

If it terminates with a Mori fibre space, then we have $n(K_{X}+B)\leq 2$. Then the semi-ampleness of $K_{X}+B$ follows from Theorem \ref{nefdim}. 

Otherwise, by Lemma \ref{Output} this $K_{X}$-MMP which is $(K_{X}+B)$-trivial terminates with a pair $(X^{\prime},B^{\prime})$ such that $K_{X^{\prime}}+(1-\varepsilon)B^{\prime}$ is nef for any sufficiently small rational $\varepsilon>0$.
Note that for any sufficiently small rational $\varepsilon>0$ we have $(X^{\prime},(1-\varepsilon)B^{\prime})$ is klt since $(X,(1-\varepsilon)B)$ is klt, and
$$\kappa(K_{X^{\prime}}+(1-\varepsilon)B^{\prime})=\kappa(K_{X^{\prime}}+B^{\prime})=\kappa(K_{X}+B)=0$$
by Lemma \ref{perturb}.
Moreover, $\mathrm{dim}\ \mathrm{Alb}(X^{\prime})\neq 0$ since $\mathrm{dim}\ \mathrm{Alb}(X)\neq 0$.
Hence by Theorem \ref{known}, $K_{X^{\prime}}+(1-\varepsilon)B^{\prime}$ is $\Q$-linearly trivial for any sufficiently small rational $\varepsilon>0$. Then $K_{X^{\prime}}+B^{\prime}$, and hence $K_{X}+B$ are $\Q$-linearly trivial.
\eproof

\btheorem\label{kltimplylc} Let $k$ be an algebraically closed field of characteristic $>3$. 
Assume we have

\noindent (1) abundance for terminal threefolds over $k$ holds, and

\noindent (2) any effective nef divisor $D$ on any klt Calabi-Yau threefold pair $(Y,\Delta)$ ($(Y,\Delta)$ is klt and $K_{Y}+\Delta\sim_{\Q}0$) over $k$ is semi-ample.

\noindent Then the abundance conjecture for threefold pairs over $k$ holds.
In particular, the abundance conjecture for klt threefold pairs over $k$ implies the abundance conjecture for lc threefold pairs over $k$.
\etheorem
\bproof Let $(X,B)$ be a projective lc threefold pair over $k$ such that $K_{X}+B$ is nef.
After replacing $(X,B)$ by its dlt modification, we can assume that $(X,B)$ is a $\Q$-factorial dlt pair and $X$ is terminal by Theorem \ref{dltmodification}.
Moreover, by Theorem \ref{nonvan} and Theorem \ref{mainthm} we can assume that $\kappa(K_{X}+B)=0$.

By Corollary \ref{MMPwithscaling} we run a $K_{X}$-MMP with scaling of $B$. It terminates by Lemma \ref{terminaltermination} since $X$ is terminal.  Hence we have a following sequence
$$(X_{0},B_{0}):=(X,B)\overset{\mu_{1}}\dashrightarrow (X_{1},B_{1})\overset{\mu_{2}}\dashrightarrow\cdots\overset{\mu_{r}}\dashrightarrow (X_{r},B_{r})$$
such that $\mu_{i}$ are $K_{X_{i-1}}$-MMP which are $(K_{X_{i-1}}+\lambda_{i-1} B_{i-1})$-trivial, where $\lambda_{i}$ are the smallest numbers such that $K_{X_{i}}+\lambda_{i}B_{i}$ are nef and $\lambda_{0}>\lambda_{1}>\cdots>\lambda_{r}$. Moreover, $(X_{r},B_{r})$ is the output of the $K_{X}$-MMP with scaling of $B$.

If $(X_{r},B_{r})$ is a minimal model, then $K_{X_{r}}$ is nef. By (1), $K_{X_{r}}$ is semi-ample. Note that
$$\kappa(X_{r},K_{X_{r}})\leq \kappa(X_{r},K_{X_{r}}+B_{r})=\kappa(X,K_{X}+B)=0.$$
Hence $K_{X_{r}}\sim_{\Q}0$. Since $K_{X_{r}}+\lambda B_{r}$ is nef for any $\lambda_{r-1}>\lambda> \lambda_{r}=0$, we have $B_{r}$ is nef on $X_{r}$. By (2), we have $B_{r}$ is semi-ample, and hence $B_{r}=0$ since 
$$\kappa(B_{r})=\kappa(K_{X_{r}}+B_{r})=\kappa(K_{X}+B)=0.$$
 It implies that $B=0$ by a standard argument using the negativity lemma (see the proof of Lemma \ref{perturb} for example).
Hence $K_{X}+B=K_{X}\sim_{\Q}0$ by (1).

Otherwise, $(X_{r},B_{r})$ is a Mori fibre space. Then we have $n(X_{r}, K_{X_{r}}+\lambda_{r}B_{r})\leq 2$, where $\lambda_{r}>0$. Hence $K_{X_{r}}+\lambda_{r}B_{r}$ is semi-ample by Theorem \ref{nefdim}. Moreover,  $K_{X_{r}}+\lambda_{r}B_{r}\sim_{\Q}0$ since 
$$\kappa(X_{r},K_{X_{r}}+\lambda_{r}B_{r})\leq \kappa(X_{r},K_{X_{r}}+B_{r})=\kappa(X,K_{X}+B)=0.$$
If $\lambda_{r}=1$, then $\lambda_{0}=\lambda_{r}=1$. It is to say that $K_{X}+B=K_{X_{r}}+\lambda_{r}B_{r}\sim_{\Q}0$. Therefore, we can assume that $\lambda_{r}<1$.
Then we have
$K_{X_{r}}+\lambda B_{r}$ is nef for any $\lambda_{r-1}>\lambda> \lambda_{r}$, and hence $B_{r}$ is nef on $X_{r}$. By (2), we have $B_{r}$ is semi-ample, and hence $B_{r}=0$ since 
$$\kappa(B_{r})=\kappa(K_{X_{r}}+B_{r})=\kappa(K_{X}+B)=0.$$
It is impossible since $\lambda_{r}>0$.
In conclusion, we have $K_{X}+B$ is semi-ample.
\eproof

\bibliographystyle{plain}
 \bibliography{Reference}

\end{document}